\documentclass[12pt]{article}

\usepackage{mathrsfs}
\usepackage{amsmath}
\usepackage{amsfonts}
\usepackage{amssymb}

\usepackage{amsmath, amsthm, amsfonts, amssymb}
\def\endpf{\hbox{\vrule height1.5ex width.5em}}

\newcommand{\s}{\mathcal}
\newcommand{\sm}{\setminus}
\newcommand{\CC}{{\mathbb C}}
\newcommand{\RR}{{\mathbb R}}
\newcommand{\BB}{{\mathbf B}}
\newcommand{\ZZ}{{\mathbb Z}}

\def \-{\bar}
\newcommand{\p}{\partial}

\newcommand{\pp}{\parallel}
\newcommand{\dbar}{\bar\partial}

\newtheorem{theorem}{Theorem}[section]
\newtheorem{lemma}[theorem]{Lemma}
\newtheorem{corollary}[theorem]{Corollary}

\newtheorem{definition}[theorem]{Definition}
\newtheorem{example}[theorem]{Example}
\newtheorem{remark}[theorem]{Remark}

\date{}

\begin{document}

\title{\bf On the existence of solutions to nonlinear systems of higher order Poisson type }

\author{Yifei Pan \  and \ Yuan Zhang\footnote{Supported in part by National Science Foundation DMS-1200652.}}

\vspace{3cm} \maketitle

\begin{abstract}
{\small In this paper, we study the existence  of higher order Poisson type systems. In detail, we prove a Residue type phenomenon for the fundamental solution of Laplacian in $\RR^n, n\ge 3$. This is analogous to the Residue theorem for the Cauchy kernel in $\CC$. With the aid of the Residue type formula for the fundamental solution, we derive the higher order derivative formula for the Newtonian potential and  obtain its appropriate $\s C^{k, \alpha}$ estimates. The  existence of  solutions to higher order Poisson type nonlinear systems is concluded as an application of the fixed point theorem.}
\end{abstract}

\bigskip
\section{Introduction  and background}

 We study  the  existence of solutions ${u}=(u_1, \ldots, u_N)$ to the following nonlinear system in $\RR^n$, $n\ge 3$:
\begin{equation}\label{system}
\triangle^m  u(x)=   a(x,   u, \nabla  u, \ldots, \nabla^{2m}   u).
\end{equation}
 Here $m\in \ZZ^+$, $\nabla^j u$ represents all $j$-th order partial derivatives of all the components of $u$, and $a:=(a_1, \ldots, a_N)$ is a vector-valued function on $x$ and the derivatives of $u$ up to order $2m$. Label the variables of $a$ by $(p_{-1}, p_0, p_1,\ldots, p_{2m})$, with $p_{-1}$ representing the position of the variable $x$ and $p_j$ representing the position of $\nabla^j  u$, $0\le j\le m$.

 The solvability  has been one of the central problems in the theory of  partial differential equations and has been explored widely since the counterexample of Hans Lewy \cite{Lw} in 1957.  See \cite{Ho1},  \cite{NT}, \cite{Mo}, \cite{Ho2}, \cite{BF},  \cite{Ln}, \cite{De} and the references therein. Unlike linear equations, there is in general no systematic theory about the solvability for nonlinear equations, not to mention nonlinear system. Most recently, \cite{Pan1} investigated the existence problem to (\ref{system}) through Cauchy Riemann operator in the case when $n=2$. In a consequent paper \cite{Pan2}, Pan studied the solvability when $m=1$ for general $n\ge 3$.

 Our main theorems are stated as follows.

\begin{theorem}\label{main2}
Let $a\in \s C^{1, \alpha}(0<\alpha <1)$. For any given appropriate initial data $\{c_j\}_{0\le j\le 2m-1}$, there exist  infinitely many solutions in  the class of $\s C^{2m,\alpha}$ to the initial value  system \begin{equation}\label{31}
\begin{split}
&\triangle^m  u(x)=   a(x,   u, \nabla  u, \ldots, \nabla^{2m-1}   u);\\
& u(0)=c_0;\\
&\nabla u(0)= c_1;\\
&\cdots\\
&\nabla^{2m-1} u(0)=c_{2m-1}
\end{split}
\end{equation}
in some small neighborhood of 0. Moreover, all those solutions are of vanishing order at most $2m$ and not radially symmetric.
\end{theorem}

Here we call $\{c_j\}_{0\le j\le 2m-1}$  {\bf an appropriate initial data} if it satisfies the appropriate symmetry conditions as derivatives of any vector-valued function. This is apparently necessary for the existence of solutions to the system (\ref{31}). A function $u$ in the class of $\s C^k$ is said to be {\bf of vanishing order $m$} ($m\le k$) at 0 if $\nabla^j u(0)=0$ for all $0\le j\le m-1$ and $ \nabla^m u(0)\ne 0$.

We point out, since the solutions obtained in Theorem \ref{main2} are of vanishing order at most $2m$, they are never trivial solutions.  Moreover, since the solutions are not radially symmetric, they are not obtained by possibly reducing the system into an ODE system with respect  only to the radial variable $r=|x|$.

Due to the flexibility of $a$, Theorem \ref{main2} can be used to construct local $m$-harmonic maps from Euclidean space to any given Riemannian manifold. The resulting image in the target manifold can be either smooth or singular, depending on the given initial data.

When $a$ is dependent also on $p_{2m}$ variable, we obtain the following existence theorem with some additional assumption on $a$.

\begin{theorem}\label{main}
If $a\in \s C^{2}$ and $a(0)=\nabla_{p_{2m}} a(0)=\nabla^2_{p_{2m}} a(0)=0$, then there exist infinitely many solutions in the class of $\s C^{2m,\alpha}$ ($0< \alpha<1$) to the system
 \begin{equation*}
\triangle^m  u(x)=   a(x,   u, \nabla  u, \ldots, \nabla^{2m}   u)
\end{equation*}
 in some small neighborhood of 0. Moreover, all those solutions are of vanishing order  $2m$ and not radially symmetric.
\end{theorem}

On the other hand, if the system (\ref{system}) is autonomous, i.e., independent of the variable $x$, then there exists solutions over large domains in the following sense.

\begin{theorem}\label{main3}
 If $a \in \s C^{2}$ and $a(0)=\nabla a(0)=0$, then for any $R> 0$, there exist infinitely  many solutions in  the class of $\s C^{2m,\alpha}$ to the autonomous system
\begin{equation*}
\triangle^m  u=   a(u, \nabla  u, \ldots, \nabla^{2m}   u)
\end{equation*}
 in $\{x\in\RR^n: |x|< R\}$. Moreover, all those solutions are of vanishing order $2m$ and are not radially symmetric.
\end{theorem}

We would like to point out, even though the autonomous system in Theorem \ref{main3} is itself translation invariant, none of  the solutions obtained there is a trivial translation of the radial solution, from the proof of Theorem \ref{main3}. On the other hand, the regularity of $a$ in Theorem \ref{main3} can be reduced to $\s C^{1, \alpha}$ if $a$ is in addition independent of $\nabla^{2m} u$ variable. This will be seen from the proof of Theorem \ref{main} and \ref{main3} in section 9 and  10. The fact will be used in some of the following examples.

  We also note the neighborhood in Theorem \ref{main2} where the solutions exist is necessarily small,  as  indicated by the following example  of Osserman.

  \begin{remark}Consider the initial value system in $n$-dimensions ($n\ge 3$):
  \begin{equation*}
  \begin{split}
    &\triangle u =|u|^{\frac{n+2}{n-2}};\\
    &u(0)= c_0;\\
    &\nabla u(0)=c_1.
    \end{split}
  \end{equation*}
 Theorem \ref{main2} applies to obtain some $\s C^{2, \alpha}$ solution over a small neighborhood of 0, say, $\{x\in \RR^n: |x|<R\}$. On the other hand, by a  result of  \cite{Os}, if the solution exists in $\{x\in \RR^n: |x|<R\}$ and $c_0>0$, then $R\le n u(0)^{-\frac{2}{n-2}}=n c_0^{-\frac{2}{n-2}}$. Consequently, $R\rightarrow 0$ as $c_0\rightarrow +\infty$. This does not contradict with Theorem \ref{main3} apparently, since the solutions constructed in Theorem \ref{main3} are of vanishing order $2m$ and hence $c_0=0$. 
 \end{remark}
\medskip

As a matter of fact, a large class of the systems fit into one or more of the three theorems. Especially, we have the following systems solvable.

\begin{example}
  For any  $p>1$ and any given $R>0$, the system
  \begin{equation*}
    \triangle^m u =\pm|u|^p
  \end{equation*}
   has infinitely many $\s C^{2m, \alpha} $ non-radial solutions over $\{x\in \RR^n: |x|<R\}$, as a consequence of Theorem \ref{main3}. Here $\alpha =\min\{1-\epsilon, p-1\}$ with $\epsilon$ any arbitrarily small positive number.  Those solutions are necessarily smooth after a standard bootstrap argument.
\end{example}

The following system has been well studied in the literature.
\begin{example} Let $H\in \s C^{3}$ and $H'(0)=0$. Consider the system
   \begin{equation*}
    \triangle u =\nabla \big(H(u)\big).
  \end{equation*}
According to Theorem \ref{main3}, for any $R>0$, the above system has infinitely many non-radial  solutions in $\s C^{2, \alpha}(\{x\in\RR^n: |x|< R\})$ for any $0<\alpha <1$.
  \end{example}
Indeed, a straight forward computation shows in the above example that $a(u, \nabla u)=\nabla \big(H(u)\big)=H'(u)\nabla u$, $\nabla_{p_0} \big(a(u, \nabla u)\big)= H''(u)\nabla u$ and $\nabla_{p_1} \big(a(u, \nabla u)\big)=H'(u)$ and hence the system satisfies $a\in \s C^2$ and  $a(0)=\nabla a(0)=0$. By Theorem \ref{main3}, for any $R>0$, there exist infinitely many solutions in the class of $\s C^{2, \alpha}(\{x\in\RR^n: |x|< R\})$ and none of them is radially symmetric.
\medskip

One similarly has the solvability  for the following $m$-th order Poisson type system.

\begin{example}
  Let $H\in \s C^{3}$ and $H'(0)=0$. Then for any $R>0$,
  \begin{equation*}
   \triangle^m  u =   \nabla \big(H(u, \nabla  u, \ldots, \nabla^{2m-2}   u)\big)
  \end{equation*}
  has infinitely many non-radial smooth solutions in $\s C^{2m, \alpha}(\{x\in\RR^n: |x|< R\})$   for any $0<\alpha <1$.
\end{example}
To see the solvability of the above example, a similar computation shows
\begin{equation*}a(u, \nabla  u, \ldots, \nabla^{2m-1}   u)=\nabla \big(H(u,\nabla  u, \ldots, \nabla^{2m-2}   u  )\big)=
\sum_{j=0}^{2m-2}\nabla_j H(u, \nabla  u, \ldots, \nabla^{2m-2} u)\nabla^{j+1} u,
 \end{equation*} where $\nabla_j H$ is the derivative of $H$ with respect to $\nabla^j u$ variable. Furthermore,\begin{equation*}
 \nabla_{p_0} \big(a(u, \nabla  u, \ldots, \nabla^{2m-1}   u)\big) = \sum_{j=0}^{2m-2}\nabla_{j}\nabla_0 H(u, \nabla  u, \ldots, \nabla^{2m-2} u)\nabla^{j+1} u \end{equation*}
  and for $k\ge 1$,
  \begin{equation*}
 \begin{split}
 \nabla_{p_k} \big(a(u, \nabla  u, \ldots, \nabla^{2m-1}   u)\big)
  =& \nabla_{p_k}\big(\sum_{j=0}^{2m-2}\nabla_j H(u, \nabla  u, \ldots, \nabla^{2m-2} u)\nabla^{j+1} u\big)\\
 =& \sum_{0\le j, k \le 2m-2}\nabla_{j}\nabla_k H(u, \nabla  u, \ldots, \nabla^{2m-2}   u)\nabla^{j+1} u \\
 & + \nabla_{k-1} H(u, \nabla  u, \ldots, \nabla^{2m-2}   u).\end{split}
  \end{equation*}Hence $a\in \s C^2$ and  $a(0)=\nabla a(0)=0$. By Theorem \ref{main3}, for any $R>0$, there exist infinitely many solutions in the class of $\s C^{2, \alpha}(\{x\in\RR^n: |x|< R\})$ and none of them is radially symmetric.
\medskip

The proof of the main theorems relies largely on a residue-type result for the fundamental solution $\Gamma(\cdot)$ of Laplacian in $\RR^n, n\ge 3$. This phenomenon is motivated by  the Residue theorem for holomorphic functions in $\CC$. In fact,  the Cauchy integral formula specially  indicates  the integral of the Cauchy kernel with a degree  $k$ holomorphic polynomial over a simple closed curve in $\CC$ is necessarily a holomorphic polynomial of the same degree. In this paper, we show the similar phenomenon also holds for $\Gamma(\cdot)$. Precisely speaking, denoting by $\s P_k$ the space of polynomials of degree $k$ restricted in  $\{x\in\RR^n: |x|< R\}$, we have, making use of zonal spherical harmonics,

\begin{theorem}\label{residue}
For any $f\in \s P_k$ with $k\ge 0$,
\begin{equation*}
  \int_{|y|=R}\Gamma(\cdot-y)f(y)d\sigma_y \in \s P_k.
\end{equation*}
Here $d\sigma_y$ is the surface area element of $\{y\in\RR^n: |y|=R\}$.
\end{theorem}

The above theorem plays an essential role in deriving the higher order derivatives for the Newtonian potential and the corresponding estimates via an induction process. Furthermore,  as a side product, the residue formula allows to define the principle value of the higher order derivatives of the Newtonian potentials.

The rest of the paper is outlined as follows. The notations  for the function spaces and the corresponding norms  are given in Section 2. In Section 3, we prove Theorem \ref{residue}. As a consequence, the principle value of higher order derivatives of the Newtonian potential is  well defined and computed in Section 4. As another application of the residue-type phenomenon for the fundamental solution, we derive the general higher order derivative formula and the corresponding estimates for the Newtonian potential in Section 5 and Section 6. Section 7 is devoted to the construction of the contraction map with the corresponding estimates necessary for the application of the fixed point theorem following the idea of \cite{Pan2}. After a delicate chasing of the parameters, we show the main theorems hold in the last three sections.

 In Appendix A, a formula of the higher order derivative of the Newtonian potential over any general domain is derived making use of the same argument as in \cite{GT}. Appendix B computes an interesting integral concerning the fundamental solution over the sphere making use of Gegenbauer polynomials. This provides a practical way to compute all the residue-type formulas for the fundamental solution.


\section{Notations}
Denote by $\BB_R$ the open ball centered at the origin with radius $R$ in $\RR^n, n\ge 3$, and denote by $\partial \BB_R$ its boundary. Namely, $\BB_R =\{x\in \RR^n: |x|<R\}$ and $\partial\BB_R =\{x\in \RR^n: |x|=R\}$. Here $|\cdot|$ is the standard Euclidean norm. We consider the following function spaces and norms over $\BB_R$ following \cite{Pan2}.

 Let $\s C(\BB_R)$ be the set of continuous functions in $\BB_R$ and $\s C^\alpha(\BB_R)$ the h\"{o}lder space in $\BB_R$ with order  $\alpha$. For $f\in \s C^\alpha(\BB_R)$, the norm of $f$ is defined by
\begin{equation*}
  \pp f\pp_{\alpha}: = \pp f\pp+R^\alpha H_\alpha[f],
\end{equation*}
where
\begin{equation*}
  \begin{split}
  &\pp f\pp: =\sup\{|f(x): x\in \BB_R\};\\
  &H_\alpha[f]: =\sup\big\{\frac{|f(x)-f(x')|}{|x-x'|^\alpha}: x, x'\in \BB_R\big\}
  \end{split}\end{equation*}
if $\pp f\pp_{\alpha}$ is finite.  We note when $\pp f\pp_{\alpha}$ is finite, the trivial extension of $f$ onto $\bar\BB_R$ is then in $\s C^\alpha(\bar\BB_R)$. $\s C^\alpha(\bar\BB_R)$ is a Banach space under the norm $\pp\cdot\pp_\alpha$.

For $k>0$, denote by $\s C^k(\BB_R)$ the collection of all functions in $\BB_R$ whose partial derivatives exist and is continuous up to order $k$. Denote by$\|\cdot\|_{\s C^k}$ the corresponding norm, where
\begin{equation*}
 \|f\|_{\s C^k}: =\sup\{\pp D^\beta f\pp: |\beta|=k\}
\end{equation*}
if $\|f\|_{\s C^k}$ is finite.

$\s C^{k, \alpha}(\BB_R)$ is  the subset of $\s C^k(\BB_R)$ whose $k$-th order derivatives  belong to $\s C^\alpha(\BB_R)$. For any multi-index $\beta=(\beta_1,\ldots, \beta_n)$ with nonnegative entries, define $|\beta|:=\sum_{j=1}^n\beta_j$ and $\beta!: = \beta_1!\cdots\beta_n!$. Given any $f\in \s C^k(\BB_R)$, we represent $D^\beta f:=\p_1^{\beta_1}\p_2^{\beta_2}\cdots\p_n^{\beta_n}f$ with $\p_j$ the partial derivative with respect to $x_j$ variables. If  $f\in\s C^{k, \alpha}(\BB_R)$, we define the semi-norm
\begin{equation*}
  \pp f\pp^{(k)}_\alpha: = \sup\{\pp D^\beta f\pp_\alpha: |\beta|=k\}
\end{equation*}
if $\pp f\pp^{(k)}_\alpha$ is finite.

Of special interest, we introduce the subset of  $\s C^{k, \alpha}(\BB_R)$ as follows.
\begin{equation*}
  \s C^{k, \alpha}_0(\BB_R): =\{f\in \s C^{k, \alpha}(\BB_R): D^\beta f(0)=0, |\beta|\le k-1\}.
  \end{equation*}

The following lemmas play an important role in the rest of the paper. The proof can be found in \cite{Pan2} and is omitted here.
\begin{lemma}\label{40}
[Pan2] If $f\in \s C^{k, \alpha}(\bar\BB_R)$, then for any $x, x'\in\bar\BB_R$ and $0< \alpha< 1$,
\begin{equation*}
  |f(x')-T_k^x(f)(x')|\le C |x-x'|^{k+\alpha}\big(\sum_{|\mu|=k}H_\alpha[D^\mu f]\big).
\end{equation*}
Here $T_k^x(f)(x')$ is the $k$-th order power series expansion of $f$ at $x$.
\end{lemma}

\begin{lemma}\label{10}
[Pan2] If $f\in \s C_0^{k, \alpha}(\bar\BB_R)$, then for any $l\le k$ and $0< \alpha< 1$,
\begin{equation*}
  \pp f\pp_\alpha^{(l)}\le C R^{k-l}\pp f\pp_\alpha^{(k)}.
\end{equation*}
\end{lemma}

\begin{remark}
  Lemma \ref{10} can be shown  to hold  for $\alpha =0$.
\end{remark}

 As a consequence of Lemma \ref{10}, $\s C^{k, \alpha}_0(\bar\BB_R)$ ($0<\alpha<1$) becomes a Banach space under the norm $\pp\cdot\pp_\alpha^{(k)}$.

Here and in the rest of the paper,  we use $C$ to represent any  positive constant number dependent only on $n, \alpha$ and $N$, where $0<\alpha<1$, $n\ge 3$ and $N\ge 1$. Especially, we point out that $C$ is independent of $R$, which is later on a key parameter in the proof of the paper.

\section{Residue-type theorem for the fundamental solution in $\RR^n$}
In complex analysis, the Residue theorem or the Cauchy integral formula states, for any holomorphic function $f$ in $\BB_R\subset\CC$ and any integer $k\ge 1$, $z\in \BB_R$,
\begin{equation*}
 \int_{|\xi|=R}\frac{f(\xi)}{\xi-z}d\xi = 2\pi i f(z),
\end{equation*}
and hence
\begin{equation*}
 \int_{|\xi|=R}D_z^k\big(\frac{1}{\xi-z}\big)f(\xi)d\xi = k!\int_{|\xi|=R}\frac{f(\xi)}{(\xi-z)^{k+1}}d\xi = 2\pi i f^{(k)}(z),
\end{equation*}
where $f^{(k)}$ is the $k$-th derivative of $f$ with respect to $z$. Recall the holomorphic kernel $\frac{1}{z}$ is the Cauchy kernel for $\dbar$ operator in $\CC$ and is also related to the first derivative of the fundamental solution in $\RR^2$. As a special case, if $f$ is a holomorphic polynomial of degree $k$ in $\BB_R\subset\CC$, then for  $z\in \BB_R$,
\begin{equation}\label{29}
  \int _{|\xi|=R}D_z^{k+1}\big(\frac{1}{\xi-z}\big)f(\xi)d\xi= 0.
\end{equation}
When $n\ge 3$, there is no holomorphic kernel for $\dbar$ operator in general.

On the other hand, the fundamental solution of Laplacian in $\RR^n, n\ge 3$ is
\begin{equation*}
  \Gamma(x)=c_n\frac{1}{|x|^{n-2}}.
\end{equation*}
Here $c_n=\frac{1}{n(2-n)\omega_n}$, with $\omega_n$ the surface area of the unit sphere in $\RR^n$. With respect to the fundamental solution, the Newtonian potential of $f$ in $\BB_R$ is defined by \begin{equation*}
  \s N(f)(x) : = \int_{\BB_R}\Gamma(x-y)f(y)dy
\end{equation*}
for any integrable function $f$ in $\BB_R$  and for $x\in \BB_R$. The Newtonian potential has attracted great attention in  physics,  and there have been many references concerning it， for instance, \cite{GT} and \cite{NW}.

The proof of  Theorem \ref{residue} makes use of zonal spherical harmonics $Z_x^{(l)}$ and its reproducing property for spherical harmonics. See \cite{SW} for the reference. In detail, let $H_l$ be the set of all spherical harmonics of degree $l$, then for any $f\in H_l$,
\begin{equation*}
 f(x)=\int_{\p\BB_1}Z_x^{(l)}(y)f(y)d\sigma_y.
\end{equation*}
Moreover, if $f\in H_k$ with $l\ne k$, then
\begin{equation*}
 0=\int_{\p\BB_1}Z_x^{(l)}(y)f(y)d\sigma_y.
\end{equation*}

  On the other hand, denote 
  by $\s P^h_k$ the space of all homogeneous polynomials of degree $k$ restricted in $\BB_R$. For any $f\in \s P^h_k$, there exist $P_j$'s, some homogenous harmonic polynomials of degree $j$, such that
  \begin{equation}\label{decompeven}
    f(x)= P_k(x) + |x|^2P_{k-2}(x) +\cdots+|x|^kP_0(x), \ \text{when $k$ is even,}
  \end{equation}
and
\begin{equation} \label{decompodd}
  f(x)= P_k(x) + |x|^2P_{k-2}(x) +\cdots+|x|^{k-1}P_1(x), \ \text{when $k$ is odd.}
\end{equation}
Note $P_j |_{\p\BB_1}\in H_j$. We now are in a position to prove the residue-type Theorem \ref{residue}  for the fundamental solution $\Gamma$ in $\RR^n$.

\medskip

{\noindent \bf Proof of Theorem \ref{residue}}: Without loss of generality, we assume $f$ is a monomial of degree $k$. We also assume that $R=1$. This is due to the following simple fact that for any $f\in \s P^h_k$,
\begin{equation*}
  \int_{\p\BB_R}\Gamma(x-y)f(y)d\sigma_y = R^{k+1}\int_{\p\BB_1}\Gamma(\frac{x}{R}-y)f(y)d\sigma_y.
\end{equation*}

Under the zonal spherical harmonics, we have, when $x\in \BB_1$,
\begin{equation*}
  \Gamma(x-y) =\sum_{l=0}^\infty C_{n,l}\frac{|x|^l}{|y|^{n+k-2}}Z^{(l)}_\frac{x}{|x|}(\frac{y}{|y|}),
\end{equation*}
where $C_{n, l}=\frac{2l+n-2}{(n-2)\omega_{n}}$. Letting $y\in \p\BB_1$, the above expression  for $x\in \BB_1$ simplifies as
\begin{equation}\label{decompgamma}
  \Gamma(x-y) =\sum_{l=0}^\infty C_{n,l}|x|^l Z^{(l)}_\frac{x}{|x|}({y}).
\end{equation}

When $k$ is odd, letting $y\in \p\BB_1$ and making use of (\ref{decompodd}), one has
\begin{equation}\label{decomphom}
  f(y)=P_k(y) + P_{k-2}(y) +\cdots+P_1(y)
\end{equation}
for some harmonic spherics $P_j\in H_j$. Therefore, combining (\ref{decompgamma}) and (\ref{decomphom}) together with the reproducing property of the zonal spherical harmonics, we have
\begin{equation*}
  \begin{split}
    \int_{\p\BB_1}\Gamma(x-y)f(y)d\sigma_y &= \int_{\p\BB_1} \bigg(\sum_{l=0}^\infty C_{n,l}|x|^l Z^{(l)}_\frac{x}{|x|}({y})\bigg)\bigg(P_k(y) + P_{k-2}(y) +\cdots+P_1(y)\bigg) d\sigma_y\\
    &= C_{n,k}|x|^kP_k(\frac{x}{|x|}) +C_{n,k-2}|x|^{k-2}P_{k-2}(\frac{x}{|x|})+\cdots +C_{n,1}|x|P_1(\frac{x}{|x|})\\
    &=  C_{n,k}P_k(x) +C_{n,k-2}P_{k-2}(x)+\cdots +C_{n,1}P_1(x)\in \s P_k.
  \end{split}
\end{equation*}

The case when $k$ is even can be treated similarly and is omitted here.
\endpf

\medskip

We remark that, despite of the constructive proof of the Residue-type formula in Theorem \ref{residue} for the fundamental solution, the integral can actually be computed directly. See Appendix B for a computation of the formula when $k=1$. The same method can practically be used for general $k>1$.

As an immediate consequence of Theorem \ref{residue}, we have the following corollary, analogous  to (\ref{29}) with respect to the Cauchy kernel in $\CC$.

\medskip

\begin{corollary}\label{residue1}
  For any $f\in \s P_k$ and any multi-index $\beta$ with $|\beta|\ge k+1$,
\begin{equation*}
  \int_{\p\BB_R}D_x^\beta\Gamma(x-y)f(y)d\sigma_y = 0
\end{equation*}
for $x\in\BB_R$.
\end{corollary}

As another consequence of Theorem \ref{residue}, we also have
\medskip

\begin{theorem}\label{2}
For any $f\in \s P_k$ and any multi-index $\beta$ with $|\beta|\ge k+2$,
\begin{equation}\label{1}
\int_{\BB_R\sm \BB_\epsilon(z)}D_x^\beta \Gamma(x-y) f(y) dy =0,
\end{equation}
when $x\in \BB_\epsilon(z)\subset \BB_R$. Here $\BB_\epsilon(z)$ is the ball centered at $z$ with radius $\epsilon$.
\end{theorem}
\medskip

{\noindent \bf Proof of Theorem \ref{2}}: Write $\beta=(\beta_1, \ldots, \beta_n)$. Without loss of generality, assume $R=1$, $\beta_1>0$ and $f$ is a monomial of degree $k$. Moreover we write $\beta' =(\beta_1-1, \ldots, \beta_n)$. Hence applying Stokes' Theorem on $D^{\beta'} \Gamma(x-y) f(y)$ over the domain $\BB_R\sm \BB_\epsilon(z)$, one has
\begin{equation}\label{3}
\begin{split}
  &\int_{\BB_1\sm \BB_\epsilon(z)}D_y^\beta \Gamma(y-x) f(y) dy\\
   =& - \int_{\BB_1\sm \BB_\epsilon(z)}D_y^{\beta'} \Gamma(y-x) \p_1f(y) dy + \int_{\p\BB_1}D_y^{\beta'} \Gamma(y-x) f(y)y_1 d\sigma_y \\
   &- \int_{\p\BB_\epsilon(z)}D_y^{\beta'} \Gamma(y-x) f(y)\frac{y_1-z_1}{|y-z|} d\sigma_y\\
   \end{split}
\end{equation}
Write $I: =\int_{\p\BB_1}D_y^{\beta'} \Gamma(y-x) f(y)y_1 d\sigma_y $ and $II: =\int_{\p\BB_\epsilon(z)}D_y^{\beta'} \Gamma(y-x) f(y)\frac{y_1-z_1}{|y-z|} d\sigma_y$. We show next that $I = II$ in $\BB_1$ and therefore
\begin{equation}\label{5}
 \int_{\BB_1\sm \BB_\epsilon(z)}D_y^\beta \Gamma(y-x) f(y) dy=- \int_{\BB_1\sm \BB_\epsilon(z)}D_y^{\beta'} \Gamma(y-x) \p_1f(y) dy.
\end{equation}

First note for $II$, after a change of coordinates by letting $y= z+\epsilon \tau$,
\begin{equation}\label{4}
  \begin{split}
    II =& \epsilon^{2-|\beta|}\int_{\p\BB_1}D_y^{\beta'} \Gamma(\frac{z-x}{\epsilon}+\tau) f(z+\epsilon \tau)\tau_1 d\sigma_\tau \\
    = & \epsilon^{2-|\beta|}\int_{\p\BB_1}D_y^{\beta'} \Gamma(\frac{z-x}{\epsilon}+\tau) \big(f(\epsilon \tau) + P_{k-1}(\tau)\big)\tau_1 d\sigma_\tau\\
     =& \epsilon^{2-|\beta|+k}\int_{\p\BB_1}D_y^{\beta'} \Gamma(\frac{z-x}{\epsilon}+\tau) f(\tau)\tau_1 d\sigma_\tau.
       \end{split}
\end{equation}
Here $P_{k-1}(\cdot)$ is some polynomial of degree $k-1$. The last identity is due to the fact that  $f$ is a monomial together with an  application of  Corollary \ref{residue1} onto $P_{k-1}(\tau)\tau_1$.

When $|\beta|\ge k+3$ and hence $|\beta'|\ge k+2$, $I$ and $II$ are both zero due to the Theorem \ref{residue} and we are done. When $|\beta|= k+2$, from (\ref{4}) we have
 \begin{equation*}
 II=\int_{\p\BB_1}D_y^{\beta'} \Gamma(\tau+ \frac{z-x}{\epsilon}) f(\tau)\tau_1 d\sigma_\tau.
  \end{equation*}
 On the other hand, $I$ by Theorem \ref{residue} is constant independent of $x
 \in \BB_1$ and therefore $I=II$ when $|\beta|\ge k+2$ and hence (\ref{5}) holds.

 Now applying the induction process on (\ref{5}), we get immediately, for $x\in \BB_1$,
 \begin{align*}
  \int_{\BB_1\sm \BB_\epsilon(z)}D_y^\beta \Gamma(y-x) f(y) dy
  = & - \int_{\BB_1\sm \BB_\epsilon(z)}D_y^{\beta'} \Gamma(y-x) \p_1f(y) dy \\
  = & \cdots \\
  = & C(f)\int_{\BB_1\sm \BB_\epsilon(z)}D_y^{\mu} \Gamma(y-x) dy \\
  = & 0.
 \end{align*}
Here $\mu$  is some multi-index with $|\mu|\ge 2$ and $C(f)$ is some constant dependent only on $f$ and $\beta$. \endpf
\medskip
\section{Principle value of higher order derivatives of  the Newtonian potential}

As a side product of Theorem \ref{2}, we show the principal value of the derivatives of the Newtonian potential exists.  Denote by $\s C^\infty_c(\RR^n)$ the set of smooth functions in $\RR^n$ with compact supports. For any $f\in \s C^\infty_c(\RR^n)$, recall the principle value of the derivatives of the Newtonian potential is defined as follows.

 \begin{definition}
 $p.v. \int_{\RR^n}D^\beta\Gamma(x-y)f(y)dy:  =\lim_{\epsilon\rightarrow 0}\int_{\RR^n\sm \BB_\epsilon(x)}D^\beta\Gamma(x-y)f(y)dy$.
 \end{definition}

\begin{theorem}\label{6}
 For any multi-index $\beta$ with $|\beta|=k\ge 2$, $f\in \s C_c^\infty$, $p.v.\int_{\RR^n}D^\beta\Gamma(x-y)f(y)dy$ exists for all $x\in \RR^n$. Moreover, if $supp f\subset \BB_R$, then
  \begin{equation}\label{7}
    p.v.\int_{\RR^n}D^\beta\Gamma(x-y)f(y)dy = \int_{\BB_R}D^\beta\Gamma(x-y)\big(f(y)-T^x_{k-2}(f)(y)\big)dy,   \end{equation}
  where $T^x_j(f)(y)$ is the $j$-th order power series expansion of $f$ at $x$. Moreover, the right hand integral of (\ref{7}) is independent of $R$.
\end{theorem}
\medskip

{\noindent \bf Proof of Theorem \ref{6}}: We first show the independence of the right hand side of (\ref{7}) of $R$. Indeed, for any two numbers $R'>R>0$ with $supp f\subset \BB_R\subset \BB_{R'}$ and when $x\in supp f$,
  \begin{equation*}
    \begin{split}
      &\int_{\BB_{R'}}D^\beta\Gamma(x-y)\big(f(y)-T^x_{k-2}(f)(y)\big)dy\\
      =& \int_{\BB_R}D^\beta\Gamma(x-y)\big(f(y)-T^x_{k-2}(f)(y)\big)dy + \int_{\BB_{R'}\sm \BB_R}D^\beta\Gamma(x-y)\big(f(y)-T^x_{k-2}(f)(y)\big)dy\\
      =& \int_{\BB_R}D^\beta\Gamma(x-y)\big(f(y)-T^x_{k-2}(f)(y)\big)dy - \int_{\BB_{R'}\sm \BB_R}D^\beta\Gamma(x-y)T^x_{k-2}(f)(y)dy
    \end{split}
  \end{equation*}
and hence by Theorem \ref{2}, when $x\in supp f\subset \BB_R\subset \BB_{R'}$,
\begin{equation*}
\int_{\BB_{R'}}D^\beta\Gamma(x-y)\big(f(y)-T^x_{k-2}(f)(y)\big)dy = \int_{\BB_R}D^\beta\Gamma(x-y)\big(f(y)-T^x_{k-2}(f)(y)\big)dy.
\end{equation*}
If $x\notin supp f$, then
\begin{equation*}
    \begin{split}
      &\int_{\BB_{R'}}D^\beta\Gamma(x-y)\big(f(y)-T^x_{k-2}(f)(y)\big)dy\\
      =&\int_{\BB_{R'}}D^\beta\Gamma(x-y)f(y)dy\\
      =&\int_{\BB_{R}}D^\beta\Gamma(x-y)f(y)dy\\
      =&\int_{\BB_{R}}D^\beta\Gamma(x-y)\big(f(y)-T^x_{k-2}(f)(y)\big)dy.
          \end{split}
  \end{equation*}

To prove (\ref{7}), we first note (\ref{7})  is trivially true if $x\notin supp f$.

When $x\in supp f$, the right hand side of (\ref{7}) is finite, due to the simple fact that when  $y\in \BB_1(x)$,
 \begin{equation*}
   |D^\beta\Gamma(x-y)|\le C|x-y|^{2-n+k}
 \end{equation*}
and \begin{equation*}
|f(y)-T^x_{k-2}(f)(y)|\le C(f)|x-y|^{k-1}
\end{equation*}
 with $C(f)$ some constant dependent  on $f$.

On the other side, making use of Theorem \ref{2}, one has
\begin{equation*}
  \int_{\BB_R\sm\BB_\epsilon(x)}D^\beta\Gamma(x-y)T^x_{k-2}(f)(y)dy=0.
\end{equation*}

As a combination of the above two facts, we obtain when $x\in supp{f}$,
\begin{align*}
  p.v. \int_{\RR^n}D^\beta\Gamma(x-y)f(y)dy =&\lim_{\epsilon\rightarrow 0}\int_{\BB_R\sm \BB_\epsilon(x)}D^\beta\Gamma(x-y)f(y)dy\\
   =& \lim_{\epsilon\rightarrow 0}\int_{\BB_R\sm \BB_\epsilon(x)}D^\beta\Gamma(x-y)\big(f(y)-T^x_{k-2}(f)(y)\big)dy\\
   =&\int_{\BB_R}D^\beta\Gamma(x-y)\big(f(y)-T^x_{k-2}(f)(y)\big)dy. \endpf
\end{align*}
\medskip

\section{ Higher order derivatives of the Newtonian potential}

It is well known that the second order derivatives of the Newtonian potential in general do not exist due to the nonintegrability of the fundamental solution after differentiation more than once. However, when $f$ is nice enough in the sense that $f\in \s C^\alpha(\BB_R)$, then $\s N(f)\in \s C^{2, \alpha}(\BB_R)$. Especially,  one has the following formula.
\begin{lemma}\label{Fr}
[Fr] Let $f\in \s C^\alpha(\BB_R)$. Then for any $x\in \BB_R$,
\begin{equation*}
  \p_{i}\p_j\s N(f)(x) =\int_{\BB_R} \p_{x_i}\p_{x_j}\Gamma(x-y)(f(y)-f(x))dy -\frac{\delta_{ij}}{n}f(x).
\end{equation*}
Moreover, for all $f\in \s C^\alpha(\BB_R)$,
\begin{equation*}
 \pp \s N(f)\pp_\alpha^{(2)}\le C\pp f\pp_\alpha,
\end{equation*}
whenever $\pp f\pp_\alpha$ is finite.
\end{lemma}

For higher order derivatives of the Newtonian potential, in the case when $n=2$, \cite{Pan1} studied it from the point of view of complex analysis. There are few references in the literature for $n>3$, though. In this section, we derive higher order derivatives of the Newtonian potential for functions in the appropriate spaces.

Throughout the rest of the paper, unless otherwise indicated, we always regard derivatives inside the integration as derivatives with respect to $y$ variables. For instance, inside an integral, $\p_1\Gamma(x-y): = \frac{\p (\Gamma(x-y))}{\p y_1}$ while $\p_{x_1}\Gamma(x-y): = \frac{\p \Gamma(x-y)}{\p x_1}$.

\begin{definition} For a given multi-index $\beta$ with $|\beta|=k+2$, $k\ge 0$,  $\s N_\beta: \s C^{k, \alpha}(\BB_R)\rightarrow \s C(\BB_R)$ is defined as follows.
 \begin{equation*} \s N_\beta(f)(x): =\int_{\BB_R} D^\beta_x\Gamma(x-y)\big(f(y)-T^x_k(f)(y)\big)dy,
\end{equation*}
for $f\in \s C^{k, \alpha}(\BB_R)$ and $x\in \BB_R$, where $T^x_k(f)(y)$ is the $k$-th order power series expansion of $f$ at $x$.\end{definition}

It is clear that  the operator $\s N_\beta$ is well defined over $\s C^{k, \alpha}(\BB_R)$.
\medskip

We next introduce the following notation for the convenience of the statement of the theorem. Given any two multi-indices $\beta=(\beta_1, \ldots, \beta_n)$ and $\mu=(\mu_1, \ldots, \mu_n)$, we say $\beta<\mu$ if $\beta_j\le\mu_j$ for each $1\le j\le n$ and $|\beta|<|\mu|$. Moreover, we define $\mu-\beta: = (\mu_1-\beta_1, \cdots, \mu_n-\beta_n)$ if $\beta<\mu$.  If in addition $|\mu|=|\beta|+1$ with $\p_iD^{\beta} = D^{\mu}$, we write $\mu -\beta= i$ .

\begin{definition}
 \begin{itemize}
  \item Given a multi-index $\beta$, we call $\{\beta^{(j)}\}_{j=1}^{k}$ a {\bf continuously increasing nesting of length $k$} for $\beta$ if $|\beta^{(j)}|=j$ for $1\le j\le k$ and $\beta^{(j)}<\beta^{(j+1)}\le \beta$ for $1\le j\le k-1$.
\item Given two multi-indices $\gamma$ and $\gamma'$, we say $\gamma'$ is {\bf the dual} of $\gamma$  with respect to $\beta$ if $D^\beta=D^\gamma D^{\gamma'}$.
\end{itemize}
\end{definition}

Making use of Theorem \ref{residue} together with Theorem \ref{33}, the following theorem gives the formula for higher order derivatives of the Newtonian potential.

\begin{theorem}\label{DN}
Let $f\in \s C^{k, \alpha}(\BB_R)$. Let $\beta$ be a multi-index with $|\beta|=k+2$ and $\{\beta^{(j)}\}$ a continuously increasing nesting of length $k+2$ for $\beta$. Then $D^\beta \s N(f)(x)$ exists for $x\in \BB_R$. Moreover,
\begin{equation}\label{44}
  D^\beta \s N(f) = \s N_\beta(f)-\sum_{j=2}^{k+2}\sum_{|\mu|=j-2} \frac{C(\beta^{(j-1)}, \mu, \beta^{(j)}-\beta^{(j-1)})}{\mu!} D^{\mu+\beta^{(j)'}}f.
\end{equation}
Here $\beta^{(j)'}$ is the dual of $\beta^{(j)}$ with respect to $\beta$, and $C(\beta^{(j-1)}, \mu, \beta^{(j)}-\beta^{(j-1)})$ is some constant dependent only on $(\beta^{(j-1)}, \mu, \beta^{(j)}-\beta^{(j-1)})$.
\end{theorem}

We point out, on the right hand side of (\ref{44}), the order $|\mu+\beta^{(j)'}|$ of the derivative of $f$ in the second term  is always equal to $k$ by definition.
\medskip

{\noindent \bf Proof of Theorem \ref{DN}}: $k=0$  is given by Lemma \ref{Fr}. When $k>0$,  for any multi-indices $\beta$ and $\mu$ with  $|\beta|= |\mu|+1$, one has by  Corollary \ref{residue1},
\begin{equation*}
\begin{split}
  \s I_{\BB_R}(\beta, \mu, j)(x): =&\int_{\p\BB_R}D^{\beta}_x\Gamma(x-y)(y-x)^{\mu}\nu_jd\sigma_y \\ =&\frac{1}{R}\int_{\p\BB_R}D^{\beta}_x\Gamma(x-y)y^{\mu}y_jd\sigma_y \\
  =&\frac{1}{R}D^{\beta}_x\int_{\p\BB_R}\Gamma(x-y)y^{\mu}y_jd\sigma_y\\
  =&R^{1+|\mu|}D^{\beta}_x\int_{\p\BB_1}\Gamma(\frac{x}{R}-y)y^{\mu}y_jd\sigma_y.\end{split}
\end{equation*}
According to Theorem \ref{residue}, $
  \int_{\p\BB_1}\Gamma(\frac{x}{R}-y)y^{\mu}y_jd\sigma_y$ is a  polynomial of degree $|\mu|+1$ in $x$ when $|x|< R$ and hence
\begin{equation*}
  \s I_{\BB_R}(\beta, \mu, j)(x)
  \equiv C(\beta, \mu, j)
\end{equation*}
with $C(\beta, \mu, j)$ some constant dependent only on $(\beta, \mu, j)$.
Therefore from Theorem \ref{33} by choosing $\Omega =\BB_R$ and $\Omega'=\BB_{R'}$ with $R'>R>0$, one obtains
\begin{equation*}
  \begin{split}
    D^\beta \s N(f)(x) =&\int_{\BB_{R'}} D^\beta_x \Gamma(x-y)\big(f(y)-T^x_k(f)(y)\big)dy\\ &-\sum_{j=2}^{k+2}D^{\beta^{(j)'}}\big(\sum_{|\mu|=j-2} \frac{D^{\mu}f(x)}{\mu!}  C(\beta^{(j-1)}, \mu, \beta^{(j)}-\beta^{(j-1)})\big)\\
    =&\int_{\BB_{R'}} D^\beta_x \Gamma(x-y)\big(f(y)-T^x_k(f)(y)\big)dy\\ &-\sum_{j=2}^{k+2}\sum_{|\mu|=j-2} \frac{C(\beta^{(j-1)}, \mu, \beta^{(j)}-\beta^{(j-1)})}{\mu!} D^{\mu+\beta^{(j)'}}f(x)
    \end{split}
  \end{equation*}
  for any $x\in \BB_R$ and any $R'>R$. We then get $\s N(f)\in \s C^{k+2}(\BB_R)$. Moreover,  for any $\beta$ with $|\beta|\le k+2$ and for any $x\in \BB_R$, after passing $R'$ to $R$,
  \begin{equation*}
  \begin{split}
    D^\beta \s N(f)(x)=&\s N_\beta(f)(x)-\sum_{j=2}^{k+2}\sum_{|\mu|=j-2} \frac{C(\beta^{(j-1)}, \mu, \beta^{(j)}-\beta^{(j-1)})}{\mu!} D^{\mu+\beta^{(j)'}}f(x). \endpf
    \end{split}
  \end{equation*}
\medskip

To simplify the notation, we define an operator $\s T_\beta$ by
 \begin{equation*}
 \s T_\beta(f)(x):=\sum_{j=2}^{k+2}\sum_{|\mu|=j-2} \frac{ C(\beta^{(j-1)}, \mu, \beta^{(j)}-\beta^{(j-1)})}{\mu!}D^{\mu+\beta^{(j)'}}f(x)
  \end{equation*}
 for any $f\in \s C^{k, \alpha}(\BB_R)$. Then $\s T_\beta: \s C^{k, \alpha}(\BB_R)\rightarrow \s C^{\alpha}(\BB_R)$ and
\begin{equation}\label{42}
\pp\s T_\beta(f)\pp_\alpha\le C\pp f\pp_\alpha^{(k)}.
\end{equation}

Under this definition, Theorem \ref{DN} can be rewritten as
\begin{equation}\label{39}
D^\beta \s N(f) = \s N_\beta(f)-\s T_\beta(f).
\end{equation}
for  any $f\in \s C^{k, \alpha}(\BB_R)$.
\bigskip

\section{H\"{o}lder norm of $D^\beta\s N(f)$}

In the derivation of the higher derivative formula of the Newtonian potential, the following operator shows up frequently and in itself is of an independent importance as well.

\begin{definition} Given multi-indices $\beta$ and $\beta'$ with $|\beta|=k+2$, $k\ge 0$ and  $D^\beta=\p_jD^{\beta'}$, a linear operator $\s {\tilde S}_\beta: \s C^{k, \alpha}(\BB_R)\rightarrow \s C(\BB_R)$ is defined as follows:
 \begin{equation*} \s {\tilde S}_\beta(f)(x): =\int_{\p\BB_R} D^{\beta'}_x\Gamma(x-y)f(y)\nu_jd\sigma_y
\end{equation*}
for $f\in \s C^{k, \alpha}(\BB_R)$ and $x\in \BB_R$. Here $d\sigma_y$ is the surface area element of $\p\BB_R$ with the unit outer normal $(\nu_1, \ldots, \nu_n)$. \end{definition}

We point out $\s {\tilde S}_\beta(f) = D^{\beta'}(\int_{\p\BB_R} \Gamma(\cdot-y)f(y)\nu_jd\sigma_y)$ is the counterpart  in $\RR^n$ of the derivatives of the Cauchy integral  for holomorphic functions in $\CC$. For the convenience of computation, we slightly modify the operator $\s {\tilde S}_\beta$ and define the following operator.

\begin{definition} Given multi-indices $\beta$ and $\beta'$ with $|\beta|=k+2$, $k\ge 0$ and  $D^\beta=\p_jD^{\beta'}$, a linear operator $\s S_\beta: \s C^{k, \alpha}(\BB_R)\rightarrow \s C(\BB_R)$  is defined as follows:
 \begin{equation*} \s {S}_\beta(f)(x): =\int_{\p\BB_R} D^{\beta'}_x\Gamma(x-y)\big(f(y)-T^x_k(f)(y)\big)\nu_jd\sigma_y
\end{equation*}
for $f\in \s C^{k, \alpha}(\BB_R)$ and $x\in \BB_R$. Here $T^x_k(f)(y)$ is the $k$-th order Taylor series expansion of $f$ at $x$, $d\sigma_y$ is the surface area element of $\p\BB_R$ with the unit outer normal $(\nu_1, \ldots, \nu_n)$. \end{definition}

Note that due to Corollary \ref{residue1} and $|\beta'|= k+1$, $\int_{\p\BB_R} D^{\beta'}\Gamma(x-y)T^x_k(f)(y)\nu_jd\sigma_y$ as a function of $x\in \BB_R$ is a constant independent of $R$. As a result of this, $\s {\tilde S}_\beta(f)$ and  $\s S(f) $ are differed by a constant only dependent on $T^x_k(f)$, and especially, independent of $R$.

We now are ready to prove the induction formula for the  derivatives of the Newtonian potential.

\begin{lemma}\label{DN2}
Let $f\in \s C^{k, \alpha}(\BB_R)$. let $\beta, \beta'$ be two multi-index with $|\beta|=k+2$ and  $D^\beta=\p_jD^{\beta'}$. We  have
\begin{equation*}
  D^\beta \s N(f) =  D^{\beta'}\s N(\p_j f) - \s S_\beta(f)-\s T_\beta(f).
\end{equation*}

\end{lemma}
{\noindent \bf Proof of Lemma \ref{DN2}}:
Making use of Stokes' Theorem in (\ref{39}) and Corollary \ref{residue1}, we have for $x\in \BB_R$,
\begin{equation*}
  \begin{split}
    D^\beta \s N(f)(x) =& \lim_{\epsilon\rightarrow 0}\int_{\BB_R-\BB_\epsilon(x)}D^\beta_x\Gamma(x-y)\big(f(y)-T^x_k(f)(y)\big)dy-\s T_\beta(f)(x)\\
    =& \lim_{\epsilon\rightarrow 0}\int_{\BB_R-\BB_\epsilon(x)}D_x^{\beta'}\Gamma(x-y)\p_j\big(f(y)-T^x_k(f)(y)\big)dy\\
    & - \lim_{\epsilon\rightarrow 0}\int_{\BB_R-\BB_\epsilon(x)}\p_j\bigg(D_x^{\beta'}\Gamma(x-y)\big(f(y)-T^x_k(f)(y)\big)\bigg)dy-\s T_\beta(f)(x)\\
    =& \lim_{\epsilon\rightarrow 0}\int_{\BB_R-\BB_\epsilon(x)}D_x^{\beta'}\Gamma(x-y)\big(\p_jf(y)-T^x_{k-1}(\p_j f)(y)\big)dy \\
    & - \int_{\p\BB_R}D_x^{\beta'}\Gamma(x-y)\big(f(y)-T^x_k(f)(y)\big)\nu_jd\sigma_y \\
    & + \lim_{\epsilon\rightarrow 0}\int_{\p\BB_\epsilon(x)}D_x^{\beta'}\Gamma(x-y)\big(f(y)-T^x_k(f)(y)\big)\nu_jd\sigma_y-\s T_\beta(f)(x)\\
    =& D^{\beta'}\s N(\p_j f)(x)- \s S_\beta(f)(x)-\s T_\beta(f)(x).
  \end{split}
\end{equation*}
 The last identity is because the third term is $O(\epsilon^{2-n-k-1+k+\alpha+n-1})=O( \epsilon^{\alpha})$. \endpf

\medskip
The following lemma shows the operator $\s S_\beta$ is a bounded operator from $\s C^{k, \alpha}(\BB_R)$ into $\s C^{\alpha}(\BB_R)$.

\begin{lemma}\label{8}
  let $\beta$ be a multi-index with $|\beta|=k+2$. The operator $\s S_\beta$ sends $\s C^{k, \alpha}(\BB_R)$ into $\s C^{\alpha}(\BB_R)$. Moreover, for any $f\in \s C^{k, \alpha}(\BB_R)$,
  \begin{equation*}
    \pp \s S_\beta(f)\pp_\alpha\le C\pp f\pp_\alpha^{(k)},
  \end{equation*}
  whenever $\pp f\pp_\alpha^{(k)}$ is finite.
  \end{lemma}

In order to prove Lemma \ref{8}, we need the following lemma.

\begin{lemma}\label{9}
For any $x\in \BB_1$, $0<\alpha<1$,
\begin{equation*}
  \int_{|y|=1}\frac{1}{|x-y|^{n-\alpha}}d\sigma_y\le C(1-|x|)^{\alpha-1}.
  \end{equation*}
\end{lemma}

   {\noindent \bf Proof of Lemma \ref{9}}: Assume $x=(r, 0, \ldots, 0)$ after rotation if necessary. One can assume in addition that $r\ge \frac{1}{2}$.

   Choose spherical coordinates  $y_1=\cos \theta_1, y_2=\sin\theta_1\cos\theta_2,\ldots, y_n=\sin\theta_1\sin\theta_2\cdots\sin\theta_{n-1}$, where $0\le \theta_1\le \pi, 0\le \theta_i\le 2\pi$ for $2\le i\le n-1$ and denote by $\p\BB_1^{n-1}$ the unit sphere in $\RR^{n-1}$, then we have,
   \begin{equation*}
     \begin{split}
       \int_{|y|=1}\frac{1}{|x-y|^{n-\alpha}}d\sigma_y =& \int_0^\pi\frac{\sin^{n-2}\theta_1}{[(\cos\theta_1-r)^2+\sin^2\theta_1]^{\frac{n-\alpha}{2}}}d\theta_1\int_{\p \BB^{n-1}_1}d\sigma_z\\
       =& C \int_0^\pi\frac{\sin^{n-2}\theta}{[1-2r\cos\theta+r^2]^{\frac{n-\alpha}{2}}}d\theta\\
       =& C \int_0^\pi\frac{\sin^{n-2}\theta}{[(1-r)^2+4r\sin^2\frac{\theta}{2}]^{\frac{n-\alpha}{2}}}d\theta\\
       \le& C(\int_0^{1-r}\frac{\sin^{n-2}\theta}{(1-r)^{n-\alpha}}d\theta +\int_{1-r}^\pi\frac{\sin^{n-2}\theta}{(2\sqrt{r}\sin\frac{\theta}{2})^{n-\alpha}}d\theta)\\
       =&A+B.
     \end{split}
   \end{equation*}
For $A$, since $\sin\theta\le \theta$ when $\theta>0$,
\begin{align*}
  A \le & \frac{C}{(1-r)^{n-\alpha}}\int_0^{1-r}\theta^{n-2}d\theta = \frac{C}{(1-r)^{n-\alpha}}(1-r)^{n-1} = C(1-r)^{\alpha-1}.
\end{align*}
For $B$, making use of $\sin\theta\ge C\theta$ when $0\le \theta\le \frac{\pi}{2}$ and the assumption $r\ge \frac{1}{2}$, we get
\begin{equation*}
  B \le C\int_{1-r}^\pi\frac{\sin^{n-2}\frac{\theta}{2}}{(\sin\frac{\theta}{2})^{n-\alpha}}d\theta =C\int_{1-r}^\pi\sin^{\alpha-2}\frac{\theta}{2}d\theta \le C\int_{1-r}^{\pi}\big(\frac{\theta}{2}\big)^{\alpha-2}d\theta \le C(1-r)^{\alpha-1}+C\le C(1-r)^{\alpha-1}.
\end{equation*}
The lemma is thus concluded. \endpf
\medskip

The following lemma in \cite{NW}, Appendix 6.2a plays an essential role in the  proof of Lemma \ref{8}.
\begin{lemma}  \label{41}
[NW] If $z$ and $z'$ are two points of the open unit disk in $\CC$, and $\gamma$ is the shorter segment of the circle through $z$ and $z'$ and orthogonal to the unit circle, then
\begin{equation*}
\int_{\gamma} \frac{|dw|}{(1-w\bar w)^{1-\alpha}} \le \frac{2}{1-\alpha}|z-z'|^{1-\alpha}
\end{equation*}
for $0<\alpha<1$.
\end{lemma}

 {\noindent \bf Proof of Lemma \ref{8}}: Write $g(y): = f(y)-T^x_k(f)(y)$.

  (i) The  estimate for  $\pp\s S_\beta(f)\pp$. Indeed, by Lemma \ref{40},
  \begin{equation*}
    \begin{split}
     |\int_{\p\BB_R} D_x^{\beta'}\Gamma(x-y)g(y)\nu_jd\sigma_y|\le & C\pp f\pp_\alpha^{(k)}R^{-\alpha} \int_{\p\BB_R}|y-x|^{2-n-k-1}|y-x|^{k+\alpha}d\sigma_y\\
     = &C\pp f\pp_\alpha^{(k)}R^{-\alpha} \int_{\p\BB_R}|y-x|^{1-n+\alpha}d\sigma_y\\
     =&C\pp f\pp_\alpha^{(k)} \int_{\p\BB_1}|y-\frac{x}{R}|^{1-n+\alpha}d\sigma_y\\
      \le &C\pp f\pp_\alpha^{(k)}.
    \end{split}
  \end{equation*}

  (ii) Given $x, x'\in \BB_R$, we estimate $|\s S_\beta(f)(x)-\s S_\beta(f)(x')|$. Assume without loss of generality that $x, x'$ lie on the plane $\{y_3=\cdots =y_n=0\}$ and write $x=Rz, x'=Rz'$ with $z, z'\in\BB_1$. Then
  \begin{equation*}
    \begin{split}
       \s S_\beta(f)(x)-\s S_\beta(f)(x')
      =& \int_{\p\BB_R} \big(D_x^{\beta'}\Gamma(x-y)-D_x^{\beta'}\Gamma(x'-y)\big)g(y)\nu_jd\sigma_y\\
      =& R^{-k}\int_{\p\BB_1} \big(D_z^{\beta'}\Gamma(z-y)-D_{z'}^{\beta'}\Gamma(z'-y)\big)g(Ry)\nu_jd\sigma_y.
    \end{split}
  \end{equation*}

Let $\gamma(t)=(\gamma_1(t), \gamma_2(t), 0, \ldots, 0): [0, 1]\rightarrow \{y_3=\cdots =y_n=0\}\cong \CC$ be a parametrization of the shorter segment of the circle through $z$ and $z'$ and orthogonal to the unit circle in $\CC$ with $\gamma(0)=z', \gamma(1)=z$. We then have
\begin{equation*}
    \begin{split}
       \s S_\beta(f)(x)-\s S_\beta(f)(x')
            =& R^{-k}\int_{\p\BB_1}\int_0^1 \frac{d}{dt} \big(D_\gamma^{\beta'}\Gamma(\gamma(t)-y)\big)dt g(Ry)\nu_jd\sigma_y\\
            =& R^{-k}\int_0^1 \sum_{k=1}^2\gamma'_k(t)dt\int_{\p\BB_1} \big(\p_{\gamma_k}D_\gamma^{\beta'}\Gamma(\gamma(t)-y)\big) g(Ry)\nu_jd\sigma_y.
    \end{split}
  \end{equation*}
Making use of Corollary \ref{residue1}, we have for any $0\le t\le 1$,
\begin{equation*}
  \begin{split}
    \int_{\p\BB_1} \big(\p_{\gamma_k}D_\gamma^{\beta'}\Gamma(\gamma(t)-y)\big) g(Ry)\nu_jd\sigma_y = \int_{\p\BB_1} \big(\p_{\gamma_k}D_\gamma^{\beta'}\Gamma(\gamma(t)-y)\big) \big(g(Ry)-T_{k}^{R\gamma(t)}(g)(Ry)\big)\nu_jd\sigma_y,
  \end{split}
\end{equation*}
where $T_{k}^{R\gamma(t)}(g)(y)$ is the $k$-th order power series expansion of $g$ at $R\gamma(t)$. Furthermore, by Lemma \ref{40},
\begin{equation*}
\begin{split}
|g(Ry)-T_{k}^{R\gamma(t)}(g)(Ry)|&\le C|Ry-R\gamma(t)|^{k+\alpha}\sum_{|\mu|=k}H_\alpha[D^\mu g]\\
&=  CR^{k+\alpha}|y-\gamma(t)|^{k+\alpha}\sum_{|\mu|=k}H_\alpha[D^\mu f].
\end{split}
\end{equation*}
Therefore,
\begin{equation*}
    \begin{split}
       &|\s S_\beta(f)(x)-\s S_\beta(f)(x')|\\
           \le& C\big(\sum_{|\mu|=k}H_\alpha[D^\mu f]\big)R^{-k}\int_0^1 \sum_{k=1}^2|\gamma'_k(t)|dt\int_{\p\BB_1} |\gamma(t)-y|^{2-n-k-2} R^{k+\alpha}|y-\gamma(t)|^{k+\alpha}d\sigma_y\\
           \le& C\big(\sum_{|\mu|=k}H_\alpha[D^\mu f]\big)R^{\alpha}\int_0^1 \sum_{k=1}^2|\gamma'_k(t)|dt\int_{\p\BB_1} |\gamma(t)-y|^{-n+\alpha} d\sigma_y.
    \end{split}
  \end{equation*}
Applying Lemma \ref{9} to $\int_{\p\BB_1} |\gamma(t)-y|^{-n+\alpha} d\sigma_y$ in the last expression, \begin{equation*}
  \begin{split}
     |\s S_\beta(f)(x)-\s S_\beta(f)(x')|
     \le &C\big(\sum_{|\mu|=k}H_\alpha[D^\mu f]\big)R^{\alpha}\int_0^1 \sum_{k=1}^2\frac{|\gamma'_k(t)|}{(1-|\gamma(t)|)^{1-\alpha}}dt\\
     \le &C\big(\sum_{|\mu|=k}H_\alpha[D^\mu f]\big)R^{\alpha}\int_0^1 \frac{|\gamma'(t)|dt}{(1-|\gamma(t)|^2)^{1-\alpha}}\\
     = &C\big(\sum_{|\mu|=k}H_\alpha[D^\mu f]\big)R^{\alpha}\int_{\gamma} \frac{|dw|}{(1-w\bar w)^{1-\alpha}}.
       \end{split}
\end{equation*}
Hence by Lemma \ref{41},
\begin{equation*}
 |\s S_\beta(f)(x)-\s S_\beta(f)(x')|\le C\big(\sum_{|\mu|=k}H_\alpha[D^\mu f]\big)R^{\alpha}|z-z'|^\alpha = C\big(\sum_{|\mu|=k}H_\alpha[D^\mu f]\big)|x-x'|^\alpha.
\end{equation*}
Namely,
\begin{equation*}
  H_\alpha[\s S_\beta(f)]\le C \big(\sum_{|\mu|=k}H_\alpha[D^\mu f]\big).
\end{equation*}
We finally have shown, combining (i) and (ii),
\begin{equation*}
  \pp \s S_\beta(f)\pp_{\alpha}\le C\pp f\pp_\alpha^{(k)}.\endpf
\end{equation*}

\begin{remark}
  In the proof of Lemma \ref{8}(ii) when carrying out the h\"older norm for $\s S_\beta$, a natural choice of $\gamma$ would usually be the segment connecting $z$ and $z$. However, the estimate in Lemma \ref{41} actually fails if $\gamma$ is chosen to be the segment instead of the geodesic as in \cite{NW}.
\end{remark}

Applying Lemma \ref{DN2} and Lemma \ref{8} inductively, one eventually obtains the following formula.

\begin{theorem}\label{DNF}
 Given a multi-index $\beta$ with $|\beta|=k+2$, let $\{\beta^{(j)}\}$ be a continuously increasing nesting for $\beta$ of length $k+2$ and $\beta^{(j)'}$ be the dual of $\beta^{(j)}$ with respect to $\beta$ for $2\le j\le k+2$. Then  for any $f\in \s C^{k, \alpha}(\BB_R)$,
\begin{equation*}
  D^\beta \s N(f) =  D^{\beta^{(2)}}\s N(D^{\beta^{(2)'}} f)- \sum_{j=3}^{k+2}\s S_{\beta^{(j)}}(D^{\beta^{(j)'}}f) - \s T_\beta(f),
  \end{equation*}
 in $\BB_R$. Moreover, for any $f\in\s C^{k, \alpha}(\BB_R)$, \begin{equation*}
 \pp \s N(f)\pp_\alpha^{(k+2)}\le C\pp f\pp_\alpha^{(k)},
 \end{equation*}
 and consequently, for any $m\in \ZZ^+$,
 \begin{equation*}
 \pp \s N^m(f)\pp_\alpha^{(k+2m)}\le C\pp f\pp_\alpha^{(k)},
 \end{equation*}
 whenever $\pp f\pp_\alpha^{(k)}$ is finite.
\end{theorem}

{\noindent \bf Proof of Lemma \ref{DNF}}: By Lemma \ref{DN2},
\begin{equation*}
  \begin{split}
    D^\beta \s N(f) &=  D^{\beta^{(k+1)}}\s N(D^{\beta^{(k+1)'}} f)- \s S_{\beta^{(k+2)}}(f) - \s T_\beta(f)\\
    &= D^{\beta^{(k)}}\s N(D^{\beta^{(k)'}} f)- \s S_{\beta^{(k+1)}}(D^{\beta^{(k+1)'}}f)- \s S_{\beta^{(k+2)}}(f) - \s T_\beta(f)\\
    &= \cdots\\
    &=  D^{\beta^{(2)}}\s N(D^{\beta^{(2)'}} f)- \sum_{j=3}^{k+2}\s S_{\beta^{(j)}}(D^{\beta^{(j)'}}f) - \s T_\beta(f).
  \end{split}
\end{equation*}
Hence from the above identity, for any $f\in\s C^{k, \alpha}(\BB_R)$ as long as $\pp f\pp_\alpha^{(k)}$ is finite,
\begin{equation*}
\begin{split}
 \pp \s N(f)\pp_\alpha^{(k+2)}:&= \sup_{|\beta|=k+2}\pp D^\beta \s N(f)\pp_\alpha\\
  &\le C\sup_{|\beta|=k+2}\big [\pp D^{\beta^{(2)}}\s N(D^{\beta^{(2)'}} f)\pp_\alpha + \sum_{j=3}^{k+2}\pp \s S_{\beta^{(j)}}(D^{\beta^{(j)'}}f)\pp_\alpha + \pp \s T_\beta(f)\pp_\alpha\big].
      \end{split}
\end{equation*}
Since $|\beta^{(j)}|=j$ and $|\beta^{(j)'}|=k+2-j$ from definition, by Lemma \ref{Fr}, Lemma \ref{8} and (\ref{42}), we get
\begin{equation*}
\begin{split}
 \pp \s N(f)\pp_\alpha^{(k+2)}&\le C\sup_{|\beta|=k+2}\big[\pp \s N(D^{\beta^{(2)'}} f)\pp^{(2)}_\alpha +\sum_{j=3}^{k+2}\pp D^{\beta^{(j)'}}f\pp_\alpha^{(j-2)}+\pp f\pp_\alpha^{(k)}\big]\\
 &\le C \sup_{|\beta|=k+2}\big[\pp D^{\beta^{(2)'}} f\pp_\alpha +\pp f\pp_\alpha^{(k+2-j+j-2)}+\pp f\pp_\alpha^{(k)}\big]\\
 &\le C\pp f\pp_\alpha^{(k)}.
  \end{split}
\end{equation*}
Finally, applying induction in the above expression,
\begin{equation*}
  \begin{split}
    \pp \s N^m(f)\pp_\alpha^{(k+2m)}& = \pp \s N(\s N^{m-1}(f))\pp_\alpha^{(k+2m)}\\
    & \le C \pp \s N^{m-1}(f)\pp_\alpha^{(k+2m-2)}\\
    & \le \cdots\\
    & \le C \pp f\pp_\alpha^{(k)}.
\endpf
  \end{split}
\end{equation*}

\section{Construction of the contraction map}
In this section, we construct a contraction map from the system (\ref{system}). Assume $a \in \s C^2$. For any  vector-valued function $f\in \s (C^{2m, \alpha}_0(\BB_R))^N$, introduce $\omega^{(1)}(f): = (\omega_1^{(1)}(f), \ldots, \omega_N^{(1)}(f))$ with
\begin{equation*}
  \omega_j^{(1)}(f)(x)=\int_{\BB_R}\Gamma(x-y)a_j(y,   f(y), \nabla  f(y), \ldots, \nabla^{2m}   f(y))dy
\end{equation*}
for $1\le j\le N$. According to Theorem \ref{DNF}, $\omega^{(1)}(f)\in \s (C^{2, \alpha}(\BB_R))^N$ and
\begin{equation*}
\pp\omega^{(1)}_j(f)\pp_\alpha^{(2)}\le C\pp a_j(\cdot, f, \ldots, \nabla^{2m}f)\pp_\alpha.
\end{equation*}
We define inductively $\omega^{(l)}(f)=(\omega_1^{(l)}(f), \ldots, \omega_N^{(l)}(f))$ for $1\le l\le m$ as follows. For each $1\le j\le N$ and $x\in \BB_R$,
\begin{equation*}
  \omega^{(l)}_j(f)(x): = \s N(\omega^{(l-1)}_j(f))(x).
\end{equation*}
Note that, in terms of the Newtonian potential,
\begin{equation*}
  \omega^{(l)}_j(f) = \s N^l\big(a_j(\cdot, f, \ldots, \nabla^{2m}f)\big).
\end{equation*}
Therefore, by Theorem \ref{DNF}, $\omega^{(l)}(f)\in \s (C^{2l, \alpha}(\BB_R))^N$ and
\begin{equation}\label{46}
\pp\omega^{(l)}_j(f)\pp_\alpha^{(2l)}\le \pp a_j(\cdot, f, \ldots, \nabla^{2m}f)\pp_\alpha.
\end{equation}

We also define $\theta(f): = (\theta_1(f), \ldots, \theta_N(f))$ from $\omega^{(m)}(f)$ by truncating degree less than $2m$ terms and part of the degree $2m$ terms in its power series expansion at 0. Precisely speaking, for $1\le j\le N$ and $x\in \BB_R$,
\begin{equation}\label{11}
  \theta_j(f)(x)=\omega_j^{(m)}(f)(x) - T_{2m-1}(\omega_j^{(m)}(f))(x) - \sum_{\beta\in \Lambda} \frac{D^{\beta}(\omega_j^{(m)}(f))(0)}{\beta!}x^\beta,
\end{equation}
where $T_{2m-1}(\omega_j^{(m)}(f))$ is the $(2m-1)$-th power series expansion of $\omega_j^{(m)}(f)$ at 0, $\Lambda=\{\beta: |\beta|=2m, \ \text{ and at least one of} \ \beta_j \ \text{is odd for} \ 1\le j\le n\}$.

From the construction, it is immediate to see that for any $f\in \s (C^{2m, \alpha}_0(\BB_R))^N$, $\triangle^m \theta(f)(x) = a(x,   f(x), \nabla  f(x), \ldots, \nabla^{2m}   f(x))$ when $x\in\BB_R$. Moreover, $\omega^{(m)}(f)\in \s (C^{2m, \alpha}(\BB_R))^N$ and so $\theta(f)\in (\s C_0^{2m,\alpha}(\BB_R))^N$. We note that, because of (\ref{46}) and (\ref{11}),  $\theta(f)$ is automatically in $(C^{2m, \alpha}_0(\bar\BB_R))^N$ after a trivial extension onto $\bar\BB_R$ if $f\in(C^{2m, \alpha}_0(\bar\BB_R))^N$.

Recall $(\s C_0^{2m, \alpha}(\bar\BB_R), \pp\cdot\pp_\alpha^{(k)})$ is a Banach space.
We now have constructed an operator between two Banach spaces as follows.
\begin{equation*}
\theta: (C^{2m, \alpha}_0(\bar\BB_R))^N \rightarrow (C^{2m, \alpha}_0(\bar\BB_R))^N
\end{equation*}
with the corresponding norm
\begin{equation*}
  \pp f\pp_\alpha^{(2m)}=\max_{1\le j\le N}\pp f_j\pp_\alpha^{(2m)}.
\end{equation*}
The  ball of radius $\gamma$ in $C^{2m, \alpha}_0(\BB_R)^N $ is denoted by
\begin{equation*}
  \s B(R, \gamma): =\{f\in C^{2m, \alpha}_0(\BB_R))^N: \pp f\pp_\alpha^{(2m)}< \gamma\}.
\end{equation*}

On the other hand, recall a function $u\in \s C^{2k}$ is called $k$-harmonic if $\triangle^{k}u=0$. Given $h = (h_1, \ldots, h_N)$ with $h_j$ any  homogeneous $m$-harmonic polynomial of  degree $2m$  and for any $f\in (\s C^{2m, \alpha}_0(\BB_R))^N$, consider
\begin{equation*}
  \theta_h(f) = h +\theta(f).
\end{equation*}
Then $\theta_h(f)\in (\s C_0^{2m,\alpha}(\BB_R))^N$, $\triangle^m \theta_h(f)(x) = \triangle^m \theta(f)(x)=a(x,   f(x), \nabla  f(x), \ldots, \nabla^{2m} f(x))$ in $\BB_R$ while  the $2m$ jets $D^\beta\theta_h(f)(0)$ with $\beta\in\Lambda$ coincide with those of the given $h$.

We will seek the solutions to (\ref{system}) by making use of the fixed point theorem. Indeed, we first show there exists $\gamma>0$ and $R>0$, such that $\theta: \s B(R, \gamma) \rightarrow \s B(R, \frac{\gamma}{2})$ and $\theta$ is a contraction map.  We then pick some nontrivial $h$ as above with $h\in \s B(R, \frac{\gamma}{2})$ and consider the corresponding operator $\theta_h$. Consequently, $\theta_h:  \s B(R, \gamma) \rightarrow \s B(R, \gamma)  $ and is a contraction map. As an application of the fixed point theorem, there exists some $u\in (\s C_0^{2m,\alpha}(\BB_R))^N$ such that $\theta_h(u)=u$. This $u$ apparently satisfies $ \triangle^m u =  \triangle^m \theta_h(u) =a(\cdot, u, \nabla u,\ldots, \nabla^{2m}u)$ in $\BB_R$ and hence is a solution to (\ref{system}) over $\BB_R$.

\begin{remark}
The construction of $\theta_h$ guarantees the solution $u$ obtained from the fixed point theorem is not a trivial solution. More precisely, $u$ is of vanishing order $2m$.
\end{remark}

We divide our proof into two steps. In each of the steps, we need to utilize Theorem \ref{DNF} and then the estimates in \cite{Pan2}.

\subsection{Estimate of $\pp \theta(f)-\theta(g)\pp_\alpha^{(2m)}$}

First, we note from (\ref{11}) that for $1\le j\le N$, for any $f, g\in \s B(R, \gamma)$,
\begin{equation*}
\begin{split}
  \pp \theta_j(f)-\theta_j(g)\pp_\alpha^{(2m)}\le& \pp \omega_j^{(m)}(f) -\omega_j^{(m)}(g)\pp_\alpha^{(2m)}+ \pp\nabla^{2m}(\omega_j^{(m)}(f)-\omega_j^{(m)}(g))\pp\\
  \le& \pp \omega_j^{(m)}(f) -\omega_j^{(m)}(g)\pp_\alpha^{(2m)}\\
  =& \pp \s N(\omega^{(m-1)}_j(f)) -\s N(\omega^{(m-1)}_j(g))\pp_\alpha^{(2m)}\\
  =& \pp \s N \big(\omega^{(m-1)}_j(f) -\omega^{(m-1)}_j(g)\big)\pp_\alpha^{(2m)}\\
  =& \cdots\\
  =& \pp \s N^{m} \big(a_j(\cdot,   f(\cdot), \nabla  f(\cdot), \ldots, \nabla^{2m} f(\cdot))-a_j(\cdot,   g(\cdot), \nabla  g(\cdot), \ldots, \nabla^{2m} g(\cdot))\big)\pp_\alpha^{(2m)}.  \end{split}
\end{equation*}

Making use of Theorem \ref{DNF} into the above expression, we have then
\begin{equation}\label{13}
\begin{split}
  \pp \theta_j(f)-\theta_j(g)\pp_\alpha^{(2m)}
  \le& C\pp a_j(\cdot,   f(\cdot), \nabla  f(\cdot), \ldots, \nabla^{2m} f(\cdot))-a_j(\cdot,   g(\cdot), \nabla  g(\cdot), \ldots, \nabla^{2m} g(\cdot))\pp_\alpha.
  \end{split}
\end{equation}

We next proceed with the same derivation of the estimate as  in section 4.1  of \cite{Pan2} without proof.
Note that due to Lemma \ref{10}, when $f\in \s B(R, \gamma)$, then $\pp\nabla^j f\pp\le CR^{2m-j}\gamma$ for $0\le j\le 2m$. Therefore as a vector-valued function of $(p_{-1}, p_0, p_1, \ldots, p_{2m})$, $a$ takes value in $E:=\{p_{-1}\in \BB_R, p_j\in \BB_{CR^{2m-j}\gamma}, 0\le j\le 2m\}$ when $u\in \s B(R, \gamma)$.

Denote by $A_j: = \sup_E |\nabla_{p_j} a|$, $Q_j: = \sup\big\{\frac{|\nabla_{p_j} a(x)-\nabla_{p_j} a(x')|}{|x-x'|^\alpha}: x, x'\in E\big\}$  and  $L_j: =\sup_E (\nabla^2_{p_jp_{2m}} a)$ with $-1\le j\le 2m$. Therefore, for $-1\le j\le 2m$,
\begin{equation}\label{27}
\begin{split}
A_j& \le C\pp \nabla_{p_j} a\pp_{\s C(E)}\le C\pp  a\pp_{\s C^{1, \alpha}(E)}, \\
 Q_j& \le C\pp  a\pp_{\s C^{1, \alpha}(E)},\\
 L_j&\le C\pp \nabla_{p_{2m}} a\pp_{\s C^1(E)}\le C\pp a\pp_{\s C^2(E)}.
 \end{split}
\end{equation}
Here $\pp  a\pp_{\s C^{1, \alpha}(E)} =\pp a\pp_{\s C^1} + \sup\big\{\frac{|\nabla a(x)-\nabla a(x')|}{|x-x'|^\alpha}: x, x'\in E\big\}$.

\begin{lemma}\label{12}[Pan2]
For any $f, g\in \s B(R, \gamma)$, if $a\in \s C^2$,
\begin{equation*}
\pp a(\cdot,   f(\cdot), \nabla  f(\cdot), \ldots, \nabla^{2m} f(\cdot))-a(\cdot,   g(\cdot), \nabla  g(\cdot), \ldots, \nabla^{2m} g(\cdot))\pp_\alpha\le \delta(R,\gamma)\pp f-g\pp_\alpha,
\end{equation*}
where
\begin{equation}\label{21}\delta(R,\gamma)=C\sum_{j=0}^{2m}R^{2m-j}\big(A_j+R^\alpha(1+R^\alpha\gamma^\alpha+\gamma)Q_j+\gamma L_j\big).\end{equation}

Moreover, if $a$ is independent of $p_{2m}$, then when $a \in \s C^{1, \alpha}$,
\begin{equation}\label{22}
 \delta(R,\gamma)=C\sum_{j=0}^{2m-1}R^{2m-j}\big(A_j+R^\alpha(1+R^\alpha\gamma^\alpha+\gamma)Q_j\big).
\end{equation}
\end{lemma}
\medskip


We then have obtained from (\ref{13}), by using Lemma \ref{12}  that
\begin{equation}\label{18}
  \pp \theta(f)-\theta(g)\pp_\alpha^{(2m)}\le \delta(R,r)\pp f-g\pp_\alpha,
\end{equation}
with $\delta(R,\gamma)$ given in (\ref{21}) or (\ref{22}).

\subsection{Estimate of $\pp \theta(f)\pp_\alpha^{(2m)}$}
Similarly, for $f\in \s B(R, \gamma)$, $1\le j\le N$,
\begin{equation}\label{15}
\begin{split}
  \pp \theta_j(f)\pp_\alpha^{(2m)}\le& \pp \omega_j^{(m)}(f) \pp_\alpha^{(2m)}+ |\nabla^{2m}(\omega_j^{(m)}(f))|\\
 \le& \pp \omega_j^{(m)}(f) \pp_\alpha^{(2m)}\\
 =& \pp \s N^{m}\big(a_j(\cdot,   f(\cdot), \nabla  f(\cdot), \ldots, \nabla^{2m} f(\cdot))    \big) \pp_\alpha^{(2m)}\\
  \le& C\pp a_j(\cdot,   f(\cdot), \nabla  f(\cdot), \ldots, \nabla^{2m} f(\cdot)) \pp_\alpha.
 \end{split}
 \end{equation}

According to the estimate in section 4.2  of \cite{Pan2},

\begin{lemma}\label{14}[Pan2]
For any $f\in \s B(R, \gamma)$, if $a\in \s C^2$,
\begin{equation*}
\pp a(\cdot,   f(\cdot), \nabla  f(\cdot), \ldots, \nabla^{2m} f(\cdot))\pp_\alpha\le \eta(R,r),
\end{equation*}
where \begin{equation}\label{17}\eta(R,\gamma) = |a(0)|+C\big(R(A_{-1} + R^\alpha\big(1+R^\alpha\gamma^\alpha +\gamma) Q_{-1}+\gamma L_{-1}\big)+\gamma\delta(R, \gamma)\big)\end{equation}
with  $\delta(R, \gamma)$ given in (\ref{21}).

Moreover, if $a$ is independent of $p_{2m}$, then when $a\in \s C^{1, \alpha}$,
\begin{equation}\label{23}
 \eta(R,\gamma) = |a(0)|+C\big(R(A_{-1} + R^\alpha\big(1+R^\alpha\gamma^\alpha +\gamma) Q_{-1}\big)+\gamma\delta(R, \gamma)\big).
\end{equation}
with $\delta(R, \gamma)$ given in (\ref{22}).
\end{lemma}
\medskip

Combining Lemma \ref{14} and (\ref{15}), we have
\begin{equation}\label{19}
  \pp \theta(f)\pp_\alpha^{(2m)}\le \eta(R,r),
\end{equation}
with $\eta(R, \gamma)$ given in (\ref{17}) or in (\ref{23}).

\section{Proof of  Theorem \ref{main}}
We now prove a slightly more general result than the main theorems following \cite{Pan2}.
\begin{theorem}\label{20}
  Let $a\in\s C^2$ and $a(0)=0$. There is a constant $\delta(<1)$ depending only on $n, N$ and $\alpha$, such that when
  \begin{equation*}
    \begin{split}
            &|\nabla_{p_{2m}} a(0)|+|\nabla^2_{p_{2m}p_{2m}} a(0)|\le \delta,
    \end{split}
  \end{equation*}
 the system (\ref{system}) has infinitely many solutions in $\s C^{2m, \alpha}$  of vanishing order $2m$ at the origin in some small neighborhood.
\end{theorem}

{\noindent \bf Proof of Theorem \ref{20}}: Our goal is to show $\theta$ sends $\s B(R, \gamma)$ into $\s B(R, \frac{\gamma}{2})$ for some positive $R$ and $\gamma$ and  is a contraction map between $\s B(R, \gamma)$. In other words, we show there exist $\gamma>0$ and $R>0$ such that for any $f, g \in\s B(R, \gamma)$,
\begin{equation*}
\pp \theta(f)-\theta(g)\pp_\alpha^{(2m)}\le c\pp f-g\pp_\alpha^{(2m)}\ \ \text{with} \ c<1
\end{equation*}
and \begin{equation*}
\pp\theta(f)\pp_\alpha^{(2m)}< \frac{\gamma}{2}.
\end{equation*}

From (\ref{18}) and (\ref{19}), it boils down to show there exist $\gamma>0$ and $R>0$ such that
\begin{equation}\label{26}
\begin{split}
  \delta(R, \gamma)&\le c<1\\
  \eta(R,\gamma)&< \frac{\gamma}{2}.
  \end{split}
\end{equation}

Denote by $\tau: = |\nabla_{p_{2m}} a(0)|+|\nabla^2_{p_{2m}p_{2m}} a(0)|$, use $\epsilon_\gamma(R)$ to represent a constant converging to 0 as $R\rightarrow 0$ for each fixed $\gamma$, and use $\epsilon(R+\gamma)$ to represent a constant converging to 0 as both $R$ and $\gamma$ go to 0. Then by continuity of $a$,
\begin{equation*}
    A_{2m}\le \tau +\epsilon(R + \gamma),\ \
    Q_{2m}\le  C\tau +\epsilon(R+ \gamma),\ \
    L_{2m}\le  \tau +\epsilon(R+ \gamma).
 \end{equation*}
 (\ref{21}) and (\ref{17}) can hence be written as
\begin{align}
      \delta(R, \gamma) &=C_a\tau(1+\gamma) + \epsilon_\gamma(R)+\epsilon(R+\gamma), \label{24} \\
    \eta(R, \gamma)& =C_a\gamma\delta(R,\gamma)+\epsilon_\gamma(R).\label{25}
  \end{align}
with $C_a$ dependent on $\pp a\pp_{\s C^2(E)}$.

First, for each $\gamma$, choose $R_0$ such that $\epsilon_\gamma(R)\le \frac{\gamma}{4}$ when $R\le R_0$ in (\ref{25}). Then we will choose  $\gamma$ and $R$ small enough so $\delta(R, \gamma)\le c: = \min\{\frac{1}{4C_a\gamma}, \frac{1}{2}\}<1$. Indeed, by choosing $\gamma(\le 1)$ and $R$ small, we can make $\epsilon_\gamma(R)+\epsilon(R+\gamma)< \frac{c}{2}$ in (\ref{24})  and hence
\begin{equation*}
  \delta(R, \gamma)< 2C_a\tau + \frac{c}{2}.
\end{equation*}
When $\tau \le \frac{c}{8C_a}$, we thus have (\ref{26}) holds.

Now recall $\Lambda=\{\beta: |\beta|=2m, \ \text{ and at least one of} \ \beta_j \ \text{is odd for} \ 1\le j\le n\}$. For $R$ and $\gamma$ chosen as above, Pick $h(x) = bx^\beta $ with $\beta\in \Lambda$, and make $b>0$ small enough such that $\pp h\pp_\alpha^{(2m)} < \frac{\gamma}{2}$ and hence $h\in \s B(R, \frac{\gamma}{2})$. Consider the operator $\theta_h(f): = h+\theta(f)$. Then $\theta_h: \s B(R, \gamma)\rightarrow \s B(R, \gamma)$ forms a contraction map from the construction. By fixed point theorem for Banach spaces, there is some $u\in \s B(R, \gamma)$ such that $\theta_h(u)=u$. $u$ thus solves the system (\ref{system}) in the class $\s C^{2m, \alpha}$ and is of vanishing order $2m$ by the construction. \endpf
\medskip

\begin{remark}\label{re}
None of the solutions constructed in the proof of Theorem \ref{20} is  radially symmetric, i.e., none of them is obtained by reducing the system (\ref{1}) possibly into an ODE system with respect to the radial variable $r=|x|$ only. Indeed, if the solution $u(x)=u(r)\in \s C_0^{2m, \alpha}$, then near 0, $u(r) = er^{2m}+o(|r|^{2m})$ for some constant $e$. In particular, $D^\beta u(0)=0$ for all $\beta\in \Lambda$. This apparently can not happen because  from the construction, $h=b x^{\beta_0}$ with  some $\beta_0\in \Lambda$ and $D^{\beta_0} u(0) = D^{\beta_0} h(0)\ne 0$.
\end{remark}

\medskip

{\noindent \bf Proof of Theorem \ref{main}}: Theorem \ref{main}  is a consequence of Theorem \ref{20} and Remark \ref{re}. \endpf

\section{Proof of Theorem \ref{main2}}
When $c_j=0, 0\le j\le 2m-1$ and $a$ is independent of $p_{2m}$, $A_{2m}, Q_{2m}$ and $L_j(-1\le j\le 2m)$ are all 0 and so (\ref{18}) and (\ref{19}) becomes
 \begin{equation*}
   \begin{split}
     \delta(R, \gamma)&\le \epsilon_\gamma(R),\\
     \eta(R, \gamma)&\le |a(0)|+\epsilon_\gamma(R).
   \end{split}
 \end{equation*}
 Here we only need $\s C^{1, \alpha}$ regularity for $a$ from the estimates (\ref{22}) and (\ref{23}). Now we choose some positive $\gamma_0$ so that $\gamma_0> 4|a(0)|$.  Consequently, we choose $R$ sufficiently small so $\epsilon_{\gamma_0}(R)\le c: =\min\{\frac{1}{2}, \frac{\gamma_0}{4}\}<1$. Hence
 \begin{equation*}
   \begin{split}
     \delta(R, \gamma_0)& \le c< 1;\\
     \eta(R, \gamma_0)&<\frac{\gamma_0}{2}.
   \end{split}
 \end{equation*}
Applying the same strategy as in the proof of Theorem \ref{20}, we can find a solution
$u\in \s B(R, \gamma_0)$ to the ODE system (\ref{31}) which is not radially symmetric.

For general given $c_{\beta}$'s with multi-indices $\beta$, we write $T_{2m-1}(x): =\sum_{j=0}^{2m-1}\frac{c_\beta}{\beta!}x^\beta$. Consider  the new system
\begin{equation*}
  \begin{split}
    &\triangle^m  \tilde u(x)=   a(x,   \tilde u+T_{2m-1}(x), \nabla  (\tilde u+T_{2m-1}(x)), \ldots, \nabla^{2m-1}   (\tilde u+T_{2m-1}(x));\\
&D^\beta\tilde u(0)=0, \ \ 0\le |\beta|\le 2m-1.\\
  \end{split}
\end{equation*}
This is a system with all the initial values equal to $0$. We then obtain some solution $\tilde u$ in the class of $\s C^{2m, \alpha}$ in some small neighborhood of 0. Then $u =\tilde u + T_{2m-1}$  solves the system (\ref{31}) in the class of $\s C^{2m, \alpha}$  in some small neighborhood of 0. Apparently, the solution obtained in this way is of vanishing order at most $2m$. Moreover, $u$ is not radially symmetric since $\tilde u$ is not. \endpf

\section{Proof of Theorem \ref{main3}}
Since $a$ is independent of $x$ and $a(0)=0$,  $A_{-1}$, $ Q_{-1}$ and $L_{-1}$  are  0 and hence in (\ref{17}),
\begin{equation*}\eta(R, \gamma)\le C_a\gamma\delta(R, \gamma).
 \end{equation*}
In order to prove Theorem \ref{main3}, we  need to show for any fixed $R>0$, there exists  some $\gamma_0>0$ such that
 \begin{equation*}
   \begin{split}
     \delta(R, \gamma_0)<1;\\
     \eta(R, \gamma_0)<\frac{\gamma_0}{2},
   \end{split}
 \end{equation*}
 which is equivalent to showing
 \begin{equation}\label{32}
   \delta(R, \gamma_0)\le c: = \min\{\frac{1}{2}, \frac{1}{2C_a}\}<1.
 \end{equation}
 Indeed, since $\nabla a(0)=0$, we have $a\in \s C^{2, 0}_0(E)$ and hence by Lemma \ref{10}, for $0\le j\le 2m$, \begin{equation*}
  A_j\le C\pp \nabla_{p_j}a\pp_{\s C^1(E)}R^{2m-j}\gamma\le C\pp a\pp_{\s C^2(E)}R^{2m-j}\gamma.
  \end{equation*}
 On the other hand, we also have by definition, for $0\le j\le 2m$,
 \begin{equation*}
 \begin{split}
  Q_j&\le C\pp \nabla_{p_j}a\pp_{\s C^1(E)}(R^{2m-j}\gamma)^{1-\alpha}\le C\pp a\pp_{\s C^2(E)}(R^{2m-j}\gamma)^{1-\alpha},\\
  L_j&\le C\pp \nabla_{p_j}a\pp_{\s C^1(E)}\le C\pp a\pp_{\s C^2(E)} .
  \end{split}
\end{equation*}
Therefore, (\ref{21}) can be written as
\begin{equation*}
  \delta(R, \gamma) = \epsilon_R(\gamma),
\end{equation*}
where $\epsilon_R(\gamma)$ represents some function converging to 0 as $\gamma$ goes to 0 for each fixed $R>0$. (\ref{32}) is thus true and the proof of Theorem \ref{main3} is  complete. \endpf

\bigskip

\appendix
\noindent{\bf \Large Appendix}
\vspace{-.5cm}
\section{ Higher order derivatives of the Newtonian potential}
We derive the higher order derivatives of the Newtonian potential following \cite{GT}. Let $\Omega\subset \RR^n$ be bounded.
\begin{definition}
Given two multi-indices   $\beta$ and $\mu$, and $j$ with $1\le j \le n$, we define for $x\in \Omega$,
\begin{equation*}
  \s I_{\Omega}(\beta, \mu, j)(x): =\int_{\p\Omega}D^{\beta}_x\Gamma(x-y)(y-x)^{\mu}\nu_jd\sigma_y, \end{equation*}
where $d\sigma_y$ is the surface area element of $\p\Omega$ with the unit outer normal $(\nu_1, \ldots, \nu_n)$.
\end{definition}

It is clear to see that $\s I_{\Omega}(\beta, \mu, j)\in \s C^\infty(\Omega)$. 

\begin{theorem}\label{33}
  Let $\beta$ be a  multi-index  with $|\beta|=k+2$. Let $\{\beta^{(j)}\}$ be a continuously increasing nesting of length $k+2$ for $\beta$ and let $\beta^{(j)'}$ be the dual of $\beta^{(j)}$ with respect to $\beta$ for $1\le j \le k+2$. Then  given  a bounded and locally $\s C^{k,\alpha}$ function $f$  in $\Omega$ and for any $x\in \Omega$,
  \begin{equation}\label{34}
  \begin{split}
    D^\beta \s N(f)(x) =&\int_{\Omega'} D^\beta_x \Gamma(x-y)\big(f(y)-T^x_k(f)(y)\big)dy\\ &-\sum_{j=2}^{k+2}D^{\beta^{(j)'}}\big(\sum_{|\mu|=j-2} \frac{D^{\mu}f(x)}{\mu!} \s I_{\Omega'}(\beta^{(j-1)}, \mu, \beta^{(j)}-\beta^{(j-1)})(x)\big).
    \end{split}
  \end{equation}
 Here $\Omega'\supset \bar\Omega$ are such that Stokes' Theorem holds on $\Omega'$ and $f$ extends to vanish on $\Omega'\sm \Omega$.
\end{theorem}

{\noindent \bf Proof of Theorem \ref{33}}: The theorem is  proved by induction on $k$. When $k=0$, the theorem reduces to the case in \cite{GT}. Assume (\ref{34}) is true for $k=k_0\ge 0$, i.e., for any $f\in \s C^{k_0,\alpha}$, and any $\beta$ with $|\beta|=k_0+2$,
\begin{equation}\label{28}
  \begin{split}
    D^\beta \s N(f)(x) =&\int_{\Omega'} D^\beta_x \Gamma(x-y)\big(f(y)-T^x_{k_0}(f)(y)\big)dy\\ &-\sum_{j=2}^{{k_0}+2}D^{\beta^{(j)'}}\big(\sum_{|\mu|=j-2} \frac{D^{\mu}f(x)}{\mu!} \s I_{\Omega'}(\beta^{(j-1)}, \mu, \beta^{(j)}-\beta^{(j-1)})(x)\big).
    \end{split}
  \end{equation}
We want to show it is true for $k=k_0+1$. Namely, for any $\beta$ with $|\beta|=k_0+3$ and $f\in \s C^{k_0+1, \alpha}$,
\begin{equation}\label{16}
  \begin{split}
    D^\beta \s N(f)(x) =&\int_{\Omega'} D^\beta_x \Gamma(x-y)\big(f(y)-T^x_{k_0+1}(f)(y)\big)dy\\ &-\sum_{j=2}^{k_0+3}D^{\beta^{(j)'}}\big(\sum_{|\mu|=j-2} \frac{D^{\mu}f(x)}{\mu!} \s I_{\Omega'}(\beta^{(j-1)}, \mu, \beta^{(j)}-\beta^{(j-1)})(x)\big).
    \end{split}
  \end{equation}
Without loss of generality, assume $D^\beta=\p_1D^{{\beta^{(k_0+2)}}}$ with $|{\beta^{(k_0+2)}}|=k_0+2$. Let
\begin{equation*}
\begin{split}
  v_\epsilon(x) = &\int_{\Omega'}D^{{\beta^{(k_0+2)}}}_x \Gamma(x-y)\eta_\epsilon(x-y)\big(f(y)-T^x_{k_0}(f)(y)\big)dy \\
  &-\sum_{j=2}^{{k_0}+2}D^{\beta^{(k_0+2)}-\beta^{(j)}}\big(\sum_{|\mu|=j-2} \frac{D^{\mu}f(x)}{\mu!} \s I_{\Omega'}(\beta^{(j-1)}, \mu, \beta^{(j)}-\beta^{(j-1)})(x)\big),
  \end{split}
\end{equation*}
where $\eta_\epsilon(x-y)=\eta(\frac{|x-y|}{\epsilon})$ with $\eta$ some smooth increasing function such that $\eta(t)=0$ when $t\le 1$ and $\eta(t)=1$ when $t\ge 2$. Here we choose $\epsilon \le \frac{dist\{\Omega'^c, \Omega\}}{2}$. When $\epsilon\rightarrow 0$, $v_\epsilon(x)\rightarrow D^{\beta^{(k_0+2)}}\s N(f)(x)$ for all $x\in \Omega$ by induction.

Now consider
\begin{equation}\label{38}
  \begin{split}
    \p_1 v_\epsilon(x) = &-\int_{\Omega'}\p_1\big(D^{{\beta^{(k_0+2)}}}_x \Gamma(x-y)\eta_\epsilon(x-y)\big)\big(f(y)-T^x_{k_0}(f)(y)\big)dy \\ &+\int_{\Omega'}D^{{\beta^{(k_0+2)}}}_x \Gamma(x-y)\eta_\epsilon(x-y)\p_{x_1}\big(f(y)-T^x_{k_0}(f)(y)\big)dy\\
    & - \p_1\big[\sum_{j=2}^{{k_0}+2}D^{\beta^{(k_0+2)}-\beta^{(j)}}\big(\sum_{|\mu|=j-2} \frac{D^{\mu}f(x)}{\mu!} \s I_{\Omega'}(\beta^{(j-1)}, \mu, \beta^{(j)}-\beta^{(j-1)})(x)\big)\big]\\
    = & A+B-\sum_{j=2}^{{k_0}+2}D^{\beta^{(j)'}}\big(\sum_{|\mu|=j-2} \frac{D^{\mu}f(x)}{\mu!} \s I_{\Omega'}(\beta^{(j-1)}, \mu, \beta^{(j)}-\beta^{(j-1)})(x)\big).
  \end{split}
\end{equation}
Here $ A: = -\int_{\Omega'}\p_1\big(D^{{\beta^{(k_0+2)}}}_x \Gamma(x-y)\eta_\epsilon(x-y)\big)\big(f(y)-T^x_{k_0}(f)(y)\big)dy $ and $B: = \int_{\Omega'}D^{{\beta^{(k_0+2)}}}_x \Gamma(x-y)\eta_\epsilon(x-y)\p_{x_1}\big(f(y)-T^x_{k_0}(f)(y)\big)dy$. We will show as $\epsilon \rightarrow 0$, for all $x\in \Omega$,
\begin{equation}\label{45}
  \begin{split}
    A+B\rightarrow  &\int_{\Omega'} D^\beta_x \Gamma(x-y)\big(f(y)-T^x_{k_0+1}(f)(y)\big)dy\\ &-\sum_{|\mu|=k_0+1} \frac{D^{\mu}f(x)}{\mu!} \s I_{\Omega'}(\beta^{(k_0+2)}, \mu, \beta^{(k_0+3)}-\beta^{(k_0+2)})(x).
  \end{split}
\end{equation}
(\ref{38}) thus gives for $x\in \Omega$,
\begin{equation*}
  \begin{split}
    \p_1 v_\epsilon(x)\rightarrow  &\int_{\Omega'} D^\beta_x \Gamma(x-y)\big(f(y)-T^x_{k_0+1}(f)(y)\big)dy\\ &-\sum_{j=2}^{k_0+3}D^{\beta^{(j)'}}\big(\sum_{|\mu|=j-2} \frac{D^{\mu}f(x)}{\mu!} \s I_{\Omega'}(\beta^{(j-1)}, \mu, \beta^{(j)}-\beta^{(j-1)})(x)\big).
  \end{split}
\end{equation*}
and hence (\ref{16}) is concluded.

For $A$,
\begin{equation*}
  \begin{split}
    A =& - \int_{\Omega'}\p_1\big(D^{{\beta^{(k_0+2)}}}_x \Gamma(x-y)\eta_\epsilon(x-y)\big)\big(f(y)-T^x_{k_0+1}(f)(y)\big)dy \\
    & - \sum_{|\mu|=k_0+1}\frac{D^\mu f(x)}{\mu!}\int_{\Omega'}\p_1\big(D^{{\beta^{(k_0+2)}}}_x \Gamma(x-y)\eta_\epsilon(x-y)\big)(y-x)^\mu dy.
  \end{split}
\end{equation*}
Applying Stokes' Theorem to the second term of the above expression, we then have
\begin{equation*}
  \begin{split}
   A = & - \int_{\Omega'}\p_1\big(D^{{\beta^{(k_0+2)}}}_x \Gamma(x-y)\eta_\epsilon(x-y)\big)\big(f(y)-T^x_{k_0+1}(f)(y)\big)dy \\
   & - \sum_{|\mu|=k_0+1}\frac{D^\mu f(x)}{\mu!}\int_{\p\Omega'}D^{{\beta^{(k_0+2)}}}_x \Gamma(x-y)\eta_\epsilon(x-y)(y-x)^\mu \nu_1d\sigma_y\\
   & + \sum_{|\mu|=k_0+1}\frac{D^\mu f(x)}{\mu!}\int_{\Omega'}D_x^{{\beta^{(k_0+2)}}} \Gamma(x-y)\eta_\epsilon(x-y)\p_1(y-x)^\mu dy.
  \end{split}
\end{equation*}
On the other hand,
\begin{equation*}
  \begin{split}
    B= & -\int_{\Omega'}D^{{\beta^{(k_0+2)}}}_x \Gamma(x-y)\eta_\epsilon(x-y)\p_{x_1}\big(T^x_{k_0}(f)(y)\big)dy.
  \end{split}
\end{equation*}
Therefore
\begin{equation*}
  \begin{split}
    A+B = & - \int_{\Omega'}\p_1\big(D^{{\beta^{(k_0+2)}}}_x \Gamma(x-y)\eta_\epsilon(x-y)\big)\big(f(y)-T^x_{k_0+1}(f)(y)\big)dy \\
   & - \sum_{|\mu|=k_0+1}\frac{D^\mu f(x)}{\mu!}\int_{\p\Omega'}D^{{\beta^{(k_0+2)}}}_x \Gamma(x-y)\eta_\epsilon(x-y)(y-x)^\mu \nu_1d\sigma_y\\
   & + \int_{\Omega'}D_x^{{\beta^{(k_0+2)}}} \Gamma(x-y)\eta_\epsilon(x-y)\big[\sum_{|\mu|=k_0+1}\frac{D^\mu f(x)}{\mu!}\p_1(y-x)^\mu -\p_{x_1}\big(T^x_{k_0}(f)(y)\big)\big] dy\\
   = & I+II+III.
    \end{split}
\end{equation*}
As $\epsilon \rightarrow 0$, for $x\in \Omega$,
\begin{equation}\label{36}
  \begin{split}
    I \rightarrow & \int_{\Omega'}D^{\beta}_x\Gamma(x-y)\big(f(y)-T^x_{k_0+1}(f)(y)\big)dy\\
    II \rightarrow & - \sum_{|\mu|=k_0+1}\frac{D^\mu f(x)}{\mu!}\int_{\p\Omega'}D^{{\beta^{(k_0+2)}}}_x \Gamma(x-y)(y-x)^\mu \nu_1d\sigma_y \\
        &= - \sum_{|\mu|=k_0+1}\frac{D^\mu f(x)}{\mu!} \s I_{\Omega'}({\beta^{(k_0+2)}}, \mu, 1)(x).
  \end{split}
\end{equation}
For $III$, notice $T_{k_0}^x(f)(y)=\sum_{|\mu|\le k_0}\frac{D^{\mu} f(x)(y-x)^\mu}{\mu!}$, so
\begin{equation*}
\begin{split}
\p_{x_1}\big(T^x_{k_0}(f)(y)\big) =& \sum_{|\mu|\le k_0}\frac{\p_1D^{\mu} f(x)(y-x)^\mu}{\mu!}+ \sum_{|\mu|\le k_0}\frac{D^{\mu} f(x)\p_{x_1}(y-x)^\mu}{\mu!}
\end{split}
\end{equation*}
On the other hand, one observes
\begin{equation*}
  \begin{split}
   \sum_{|\mu|\le k_0}\frac{\p_1D^{\mu} f(x)(y-x)^\mu}{\mu!}=& \sum_{|\mu|=k_0+1}\frac{D^\mu f(x)}{\mu!}\p_1(y-x)^\mu
+\sum_{|\mu|\le k_0-1}\frac{\p_1D^{\mu} f(x)(y-x)^\mu}{\mu!}, \\
\sum_{|\mu|\le k_0}\frac{D^{\mu} f(x)\p_{x_1}(y-x)^\mu}{\mu!}
=&-\sum_{|\mu|\le k_0-1}\frac{\p_1D^{\mu} f(x)(y-x)^\mu}{\mu!}.
  \end{split}
\end{equation*}
Hence
\begin{equation*}
  \begin{split}
    \sum_{|\mu|=k_0+1}\frac{D^\mu f(x)}{\mu!}\p_1(y-x)^\mu -\p_{x_1}\big(T^x_{k_0}(f)(y)\big) = 0,
  \end{split}
\end{equation*}
and \begin{equation}\label{37}
  III=0.
\end{equation}
Combining (\ref{36}) and (\ref{37}), (\ref{45}) thus holds. \endpf

\section{Computation of $\s I_{B_1}(0, 0, 1)$}
We will compute $\s I_{B_1}(0, 0, 1)(x): = \int_{\p\BB_1}\Gamma(x-y)\nu_1d\sigma_y$ for $x\in \BB_1$.

Write $x = U\cdot [a, 0, \ldots, 0]^t$, where $U=(u_{ij})_{1\le i, j\le n}$ is some unitary matrix and $a=|x|$, and then make a change of coordinates by letting $y= U\cdot\tilde y$ in the expression of $\s I_{B_1}(0, 0, 1)$. We then get
\begin{equation*}
  \s I_{B_1}(0, 0, 1): = \int_{\p\BB_1}\frac{\sum_{0\le j\le n}u_{1j}\tilde y_j}{\sqrt{(a-\tilde y_1)^2+\tilde y_2^2+\cdots +\tilde y_n^2}^{n-2}}d\sigma_{\tilde y}
\end{equation*}
Write $\tilde y$ back by $y$ allowing  an abuse of notation, then
\begin{equation*}
\begin{split}
  \s I_{B_1}(0, 0, 1): = &u_{11}\int_{\p\BB_1}\frac{ y_1}{\sqrt{(a-  y_1)^2+  y_2^2+\cdots +  y_n^2}^{n-2}}d\sigma_{y}\\
  & +\sum_{2\le j\le n}u_{1j}\int_{\p\BB_1}\frac{ y_j}{\sqrt{(a- y_1)^2+  y_2^2+\cdots +  y_n^2}^{n-2}}d\sigma_{y}
  \end{split}
\end{equation*}
Since $\frac{ y_j}{\sqrt{(a- y_1)^2+  y_2^2+\cdots +  y_n^2}^{n-2}}$ is odd with respect to $y_j$ when $j\ge 2$,
$$\int_{\p\BB_1}\frac{ y_j}{\sqrt{(a- y_1)^2+  y_2^2+\cdots +  y_n^2}^{n-2}}d\sigma_{y}= 0$$ when $j\ge 2$ and hence
\begin{equation*}
 \s I_{B_1}(0, 0, 1): = u_{11}\int_{\p\BB_1}\frac{ y_1}{\sqrt{(a- y_1)^2+  y_2^2+\cdots +  y_n^2}^{n-2}}d\sigma_{y}.
\end{equation*}
Next, rewrite the above integral in terms of the spherical coordinates, then we obtain
\begin{equation}\label{43}
\begin{split}
  \s I_{B_1}(0, 0, 1) =&\omega_{n-1}u_{11}\int_{-\frac{\pi}{2}}^{\frac{\pi}{2}}\frac{\sin t\cos^{n-2} t}{\sqrt{(a-\sin t)^2+\cos^2 t}^{n-2}}dt\\
  =&\omega_{n-1}u_{11}\int_{-1}^1\frac{u(1-u^2)^{\frac{n-3}{2}}}{(1-2au+a^2)^{\frac{n-2}{2}}}du
  \end{split}
\end{equation}
Here $\omega_{n-1}$ is the surface area of the unit sphere in $\RR^{n-1}$.

In order to compute (\ref{43}), we need to use Gegenbauer polynomials. Recall for each fixed $\rho$, the Gegenbauer polynomials are $\{C_n^{(\rho)}(x)\}$ in $[-1, 1]\subset\RR$ satisfying
\begin{equation*}
  \frac{1}{(1-2xt +t^2)^\rho}=\sum_{n=0}^\infty C_n^{(\rho)}(x)t^n
\end{equation*}
in $(-1, 1)$. In particular,
\begin{equation*}
  \begin{split}
    C_0^{(\rho)}(x)&=1,\\
    C_1^{(\rho)}(x)&=2\rho x,\\
    C_n^{(\rho)}(x)&=\frac{1}{n}[2x(n+\rho-1)C_{n-1}^{(\rho)}(x)-(n+2\rho-2)C_{n-2}^{(\rho)}(x)].
  \end{split}
\end{equation*}
Moreover, $\{C_n^{(\rho)}(x)\}$ are orthogonal polynomials on the interval [-1,1] with respect to the weight function $(1 - x^2)^{\rho-\frac{1}{2}}$. In other words,
\begin{equation*}
  \begin{split}
    &\int_{-1}^1C_n^{(\rho)}(x)C_m^{(\rho)}(x)(1 - x^2)^{\rho-\frac{1}{2}}dx=0, m\ne n\\
    &\int_{-1}^1[C_n^{(\rho)}(x)]^2(1 - x^2)^{\rho-\frac{1}{2}}dx=\frac{\pi 2^{1-2\rho}\Gamma(n+2\rho)}{n!(n+\rho)\Gamma(\rho)^2}.
  \end{split}
\end{equation*}

Letting $\rho =\frac{n-2}{2}$, then
$\frac{1}{(1-2au+a^2)^{\frac{n-2}{2}}}=\sum_{n=0}^\infty C_n^{(\rho)}(u)a^n$ and $u=\frac{C_1^{(\rho)}(u)}{2\rho}$. (\ref{43}) can hence be written as
\begin{equation*}
  \begin{split}
    \s I_{B_1}(0, 0, 1) =&\omega_{n-1}u_{11}\int_{-1}^1\sum_{n=0}^\infty C_n^{(\rho)}(u)a^n \cdot\frac{C_1^{(\rho)}(u)}{2\rho}(1-u^2)^{\rho-\frac{1}{2}}du\\
    =&\omega_{n-1}u_{11}\int_{-1}^1C_1^{(\rho)}(u)a \cdot \frac{C_1^{(\rho)}(u)}{2\rho}(1-u^2)^{\rho-\frac{1}{2}}du\\
    =&\omega_{n-1}u_{11}\frac{a}{2\rho}\frac{\pi 2^{1-2\rho}\Gamma(1+2\rho)}{(1+\rho)\Gamma(\rho)^2}\\
    =&\frac{4\pi^{\frac{n}{2}}}{n\Gamma(\frac{n-2}{2})}x_1.
  \end{split}
\end{equation*}

\begin{remark}
Making use of the same approach as the above in addition to observing some symmetry of the integrand over the sphere, one can practically compute  $\s I_{\BB_R}(\beta, \mu, j)$ for all $(\beta, \mu, j)$.
\end{remark}



\end{document}